\documentclass{amsart}

\usepackage{amsthm}
\usepackage{amssymb}
\usepackage{amsmath}
\newtheorem {Theorem}   {Theorem}

\numberwithin{Theorem}{section}

\newtheorem {lemma}[Theorem]{Lemma}         
\newtheorem {proposition}[Theorem]{Proposition}  

\newtheorem{definition}[Theorem]{Definition}
\theoremstyle{remark}

\newtheorem{theorem}{Theorem}

\newtheorem{corollary}[theorem]{Corollary}

\begin{document}
\author{Alex Clark}
\title{Linear Flows on $\kappa $-Solenoids}
\address{Department of Mathematics, University of North Texas, 
Denton, TX 76203-5118}
\email{alexc@unt.edu}
\date{13 January 1998}

\begin{abstract}
Linear flows on inverse limits of tori are defined and it is shown that two
linear flows on an inverse limit of tori are equivalent if and only if there
is an automorphism of the inverse limit generating the equivalence.
\vspace{.2in}
\end{abstract}

\maketitle

\section{Introduction}

We define the families of linear flows on the inverse limits of
finite-dimensional tori $\mathbf{T}^n$ ($n$ fixed) with epimorphic bonding
maps and on a special class of inverse limits of $\mathbf{T}^\infty $ and
prove that two flows from such a family are topologically equivalent if and
only if there is an automorphism generating the equivalence. This result
generalizes the well known classification of the linear flows on $\mathbf{T}
^2$ (see \cite{I}, pp. 36-38). While there does not seem to be a proof in the
literature of this result for the general linear flow on $\mathbf{T}^\kappa $
for $\kappa >2$, some related questions are addressed in (\cite{KH}, 2.3). This
result reduces the problem of the classification of these linear flows to
the classification of the automorphisms of the corresponding inverse limit.
We give a characterization of the automorphisms on the finite product of
one-dimensional solenoids and we work out the 2-dimensional case in detail,
thereby classifying the linear flows on such products. We also find a
condition on the character group of a finite-dimensional inverse limit as
above that determines when the inverse limit is isomorphic with a product of
one-dimensional solenoids.

It can be shown that the subgroup of the reals generated by the Bohr-Fourier
exponents of an almost periodic orbit of a flow in a complete metric space
determines the equivalence class of the flow obtained by the restriction of
the original flow to the compact minimal set which the closure of the image
of the orbit forms in the sense that any two such orbits with the same
associated group determine equivalent flows (see, for example, \cite{LZ}, 3\S 2).
This, together with a straightforward application of Pontryagin duality, can
be used to demonstrate that a flow in a complete metric space restricted to
the closure of the image of an almost periodic orbit is equivalent to an
irrational linear flow as defined here. Any two (topologically) equivalent
almost periodic flows are equivalent to members of the same family of
irrational linear flows, and so our results serve as a program for the
classification of almost periodic flows in complete metric spaces. It then
follows from a theorem of Nemytskii (\cite{NS}; V, 8.16) that any metric compact
connected abelian group is isomorphic with an inverse limit of the type
treated here (a $\kappa -$solenoid) (see also \cite{Pont}, Thm 68). 
\footnote{
A more detailed exposition of these matters appears in the dissertation of
the author.}

\section{$\kappa -$Solenoids}

\subsection{$n-$Solenoids}

$S^1\overset{def}{=}\mathbf{T}^1\overset{def}{=}\mathbb{R}/\mathbb{Z}$, with group
operation $``+"$ inherited from the covering homomorphism $p^1:$ $\mathbb{
R\rightarrow }S^1;$\ $x\mapsto x\left( \text{mod}1\right) $. Unless
otherwise stated, $n$ denotes a member of the set of natural numbers $\mathbb{
N=\{}1,2,...\}$ and $\infty $ denotes the countably infinite cardinal. The $
n $-torus is denoted $\mathbf{T}^n\overset{def}{=}\prod_{i=1}^nS^1$ and $
\mathbf{T}^\infty \,$ is defined to be the space $\prod_{i=1}^\infty S^1$
and we let $\mathbf{x=}\langle x_1,...,x_i,...\rangle $ denote a point of $
\mathbf{T}^\kappa $ for $\kappa \in \mathbb{N}$ or $\kappa =\infty $, and we
give $\mathbf{T}^\kappa $ the metric $d_\kappa ;$ 
\begin{equation*}
d_\kappa \left( \mathbf{x},\mathbf{y}\right) \overset{def}{=}
\sum_{i=1}^\kappa \frac 1{2^i}\left| u_i-v_i\right| \text{,}
\end{equation*}
where $u_i,v_i\in \mathbb{R}\,$are representatives of the classes of $
x_i,y_i\,\;$chosen so that$\;\left| u_i-v_i\right| \leq \frac 12$. An
inverse limit of $\mathbf{T}^\kappa $ ($\kappa \in \mathbb{N}$ or $\kappa
=\infty $, fixed) with epimorphic bonding maps has the group structure
inherited from the Cartesian product $\prod_{i=1}^\infty \mathbf{T}^\kappa $
. For a fixed $\kappa $, $\prod_{i=1}^\infty \mathbf{T}^\kappa $ and its
subspaces are given the metric $d_\kappa ^\infty ;$ 
\begin{equation*}
d_\kappa ^\infty \left( \langle \mathbf{x}^j\rangle _{j=1}^\infty ,\langle 
\mathbf{y}^j\rangle _{j=1}^\infty \right) \overset{def}{=}\sum_{j=1}^\infty
\frac 1{2^j}d_\kappa \left( \mathbf{x}^j,\mathbf{y}^j\right) .
\end{equation*}

We assume throughout that all homomorphisms (automorphisms, etc.) between
topological groups are continuous.

\begin{definition}
For $\kappa \in \mathbb{N\,\cup }\,\mathbb{\,}\left\{ \infty \right\} $, $p^\kappa
:\mathbb{R}^\kappa \rightarrow \mathbf{T}^\kappa $ is the homomorphism\\$
\mathbf{t}=\left( t_1,...,t_i,...\right) \mapsto \langle p^1\left( \
t_1\right) ,...,p^1\left( \ t_i\right) ,...\rangle $.
\end{definition}

Notice that if $f:\mathbf{T}^n\rightarrow \mathbf{T}^n$ is a homomorphism
there is a unique homomorphism represented by a matrix with integer entries $
M:\mathbb{R}^n\rightarrow \mathbb{R}^n$ satisfying $f\circ p^n\left( \mathbf{t}
\right) =p^n\circ M\left( \mathbf{t}\right) $, and when $f$ is an
epimorphism $\det M\in \mathbb{Z}-\{0\}$. And if $f:\mathbf{T}^n\rightarrow 
\mathbf{T}^n$ is represented by the $n\times n$ integer matrix $M\,$ with $
\det M\neq 0$, there are two $n\times n$ integer matrices $P\,$ and $Q\,$
which have inverses with integer entries and which satisfy $M=P\Delta Q$,
where $\Delta $ is a diagonal matrix with integer entries. Then with $
k=\left| \det M\right| =\left| \det \Delta \right| $, $f\,$ is a $k-$to$-$
one covering map since $\Delta $ represents such a map.

\begin{definition}
For a fixed $n$ and a sequence $\overline{M}=\left( M_1,M_2,...\right) $ of $
n\times n$ matrices $M_i$ with integer entries and non-zero determinants, we
define the topological group $\sum\nolimits_{\overline{M}}$ with identity $
e_{\overline{M}}$ to be the inverse limit of the inverse sequence $\{\mathbf{
X}_j,f_j^i\}$, where $\mathbf{X}_j=\mathbf{T}^n$ for all\ $j\in \mathbb{N}$ and 
$f_j^{j+1}$ is the topological epimorphism represented by the matrix $M_j$;\ 
$f_j^{j+1}\circ p^n=p^n\circ M_j$. 
\begin{equation*}
\sum\nolimits_{\overline{M}}\overset{def}{=}\lim\limits_{\leftharpoondown }\{
\mathbf{X}_j,\text{ }f_i^j\}\subset \prod_{j=1}^\infty \mathbf{T}^n\text{, }
\end{equation*}
and we define such an inverse limit $\sum\nolimits_{\overline{M}}$ to be an 
\underline{$n$-solenoid}.
\end{definition}

\subsection{$\infty -Solenoids$}

\begin{definition}
If $f:\mathbf{T}^n\rightarrow \mathbf{T}^n;$\ $\langle x_1,...,x_n\rangle
\mapsto \langle y_1,...,y_n\rangle $ is a homomorphism represented by a
matrix with integer entries $M:\mathbb{R}^n\rightarrow \mathbb{R}^n;$\ $\left(
t_1,...,t_n\right) \mapsto \left( s_1,...,s_n\right) $, then we define the
following maps: $f\times id:\mathbf{T}^\infty \rightarrow \mathbf{T}^\infty $
; 
\begin{equation*}
\langle x_1,x_2,...,\rangle \mapsto \langle
y_1,...,y_n,x_{n+1},x_{n+2},...\rangle
\end{equation*}
and $M\times id:\mathbb{R}^\infty \rightarrow \mathbb{R}^\infty $ ; 
\begin{equation*}
\left( t_1,t_2,...\right) \mapsto \left(
s_1,...,s_n,t_{n+1},t_{n+2},...\right) .
\end{equation*}
And we define a map such as $f\times id$ to be an \underline{$n$-map}.
\end{definition}

Notice that if $M$ represents $f$, then $M\times id$ represents $f\times id$
in the sense that $\left( f\times id\right) \circ p^\infty \left( \mathbf{t}
\right) =p^\infty \circ \left( M\times id\right) \left( \mathbf{t}\right) $.
Also, if $f\times id$ \ and $g\times id$\ are both $n$-maps, then $\left(
f\times id\right) \circ \left( g\times id\right) $ is the $n$-map $\left(
f\circ g\right) \times id$. And if $f\times id$ is an $n$-map $;$ 
\begin{equation*}
\langle x_1,x_2,...,\rangle \mapsto \langle
y_1,...,y_n,x_{n+1},x_{n+2},...\rangle
\end{equation*}
and $\nu >n$, it is possible to represent $f\times id$ as the $\nu $-map $
f^{\prime }\times id$, where $f^{\prime }:\mathbf{T}^\nu \rightarrow \mathbf{
T}^\nu ;$
\begin{equation*}
\langle x_1,...,x_\nu \rangle \mapsto \langle y_1,...,y_n,x_{n+1},...,x_\nu
\rangle .
\end{equation*}
Also, if $f\times id$ is an $n$-map and $g\times id$ is an $m$-map, then $
\left( f\times id\right) \circ \left( g\times id\right) $ is a $\max \{m,n\}$
-map, for we may represent both $f\times id$ and $g\times id$ as $\max
\{m,n\}$-maps and then the above observation on compositions applies.

\begin{definition}
If for each $i\in \mathbb{N\;}$ $g_i^{i+1}\times id$ is an $n_i$-map
represented by the map $M_i\times id:\mathbb{R}^\infty \rightarrow \mathbb{R}
^\infty $, where $M_i$ is an $n_i\times n_i$ integer matrix with non-zero
determinant, we define the topological group $\sum_{\overline{M}}$ with
identity $e_{\overline{M}}$ to be the inverse limit of the inverse sequence $
\{G_i,f_i^j\}$, where $G_i=\mathbf{T}^\infty $ for all $i\in \mathbb{N}$ and $
f_i^{i+1}=$ $g_i^{i+1}\times id$; 
\begin{equation*}
\sum\nolimits_{\overline{M}}\overset{def}{=}\lim\limits_{\leftharpoondown
}\{G_i,f_i^j\}\subset \prod_{j=1}^\infty \mathbf{T}^\infty \text{,}
\end{equation*}
and we define such an inverse limit $\sum\nolimits_{\overline{M}}$ to be an 
\underline{$\infty -$solenoid}.
\end{definition}

In the following $\sum\nolimits_{\overline{M}}$ represents a $\kappa -$
solenoid for some $\kappa \in \mathbb{N}\cup \left\{ \infty \right\} $.

\begin{definition}
$f_i:\sum\nolimits_{\overline{M}}\rightarrow \mathbf{T}^\kappa ;$\ $\left( 
\mathbf{x}^1,\mathbf{x}^2,...\right) \mapsto \mathbf{x}^i$ is projection
onto the $i^{th}$ factor.
\end{definition}

\begin{definition}
$\mathbf{C}_{\overline{M}}$ is the path component of $e_{\overline{M}}$ in $
\sum\nolimits_{\overline{M}}$.
\end{definition}

\begin{definition}
\label{Pi def}If $\kappa =n<\infty $, we define $\pi _{\overline{M}\ }:\mathbb{R
}^n\rightarrow \sum_{\overline{M}}$ to be the homomorphism 
\begin{equation*}
\mathbf{t=}\left( t_1,...,t_n\right) \mapsto \left( p^n\left( \mathbf{t}
\right) ,p^n\circ M_1^{-1}\left( \mathbf{t}\right) ,...,p^n\circ
M_j^{-1}\circ \cdots \circ M_1^{-1}\left( \mathbf{t}\right) ,...\right) 
\text{.}
\end{equation*}
and if $\kappa =\infty $, we define $\pi _{\overline{M}\ }:\mathbb{R}^\infty
\rightarrow \sum_{\overline{M}}$ to be the homomorphism 
\begin{equation*}
\mathbf{t\mapsto }\left( p^\infty \left( \mathbf{t}\right) ,p^\infty \circ
\left( M_1^{-1}\times id\right) \left( \mathbf{t}\right) ,...,p^\infty \circ
\left( M_j^{-1}\times id\right) \circ \cdots \circ \left( M_1^{-1}\times
id\right) \left( \mathbf{t}\right) ,...\right) \text{,}
\end{equation*}
\end{definition}

Notice that $\ker \pi _{\overline{M}\ }\subset \ker \left( f_1\circ \pi _{
\overline{M}\ }\right) =\mathbb{Z}^\kappa $.

That $\pi _{\overline{M}\ }$ maps $\mathbb{R}^n$ onto $\mathbf{C}_{\overline{M}
} $ follows from Theorem 5.8 in [McC], where he characterizes the path
components of the inverse limit of an inverse sequence with all bonding maps
regular covering maps between spaces which admit a universal covering.
However, we shall demonstrate directly that $\pi _{\overline{M}\ }\left( 
\mathbb{R}^n\right) =\mathbf{C}_{\overline{M}}$ and generalize this result in 
\textbf{Corollary \ref{Pi}}, and to do so we prove two preliminary lemmas
using the terminology and results in (\cite{S}, Chapt 2).

\begin{lemma}
\label{Fib}Let $p_1:E\rightarrow L$ and $p_2:L\rightarrow B$ be maps which
satisfy the following conditions: $\left( \mathbf{1}\right) $ $\,$ $p_2$ has
unique path lifting $\;$and $\;\left( \mathbf{2}\right) $ $p=p_2\circ p_1$
is a fibration. Then $p_1$ is a fibration.
\end{lemma}

\textbf{Proof}: If $F,f^{\prime }$ are maps as in the following diagram, we
need to find a map $F^{\prime }:X\times \left[ 0,1\right] \rightarrow E$
(represented by the diagonal arrow in the diagram) which makes the following
diagram (A) commute:

\begin{equation*}
\begin{array}{ccc}
X\times \{0\} & \overset{f^{\prime }}{\rightarrow } & E \\ 
\downarrow ^{\cap } & \nearrow & \downarrow ^{p_1} \\ 
X\times \left[ 0,1\right] & \underset{F}{\rightarrow } & L
\end{array}.
\end{equation*}
Then for any $x\in X$, $p_1\circ f^{\prime }\left( x,0\right) =F\left(
x,0\right) $. Since $p$ is a fibration, there is a map $G:X\times \left[
0,1\right] \rightarrow E$ making the following diagram (B) commute:

\begin{equation*}
\begin{array}{ccc}
X\times \{0\} & \overset{f^{\prime }}{\rightarrow } & E \\ 
\downarrow ^{\cap } & \nearrow _G & \downarrow ^p \\ 
X\times \left[ 0,1\right] & \underset{p_2\circ F}{\rightarrow } & B
\end{array}.
\end{equation*}
Then for any $x\in X$, $p_1\circ G\left( x,0\right) =p_1\circ f^{\prime
}\left( x,0\right) =F\left( x,0\right) $ ; the first equality follows from
diagram (B) and the second equality from the observation after diagram (A).
Now fix $x\in X$ and define the paths $\omega $ and $\omega ^{\prime }$ in $
L $ : 
\begin{equation*}
\omega \left( t\right) \overset{def}{=}F\left( x,t\right) \text{ and }\omega
^{\prime }\left( t\right) \overset{def}{=}p_1\circ G\left( x,t\right) \text{.
}
\end{equation*}
Then we have $\omega \left( 0\right) =F\left( x,0\right) =p_1\circ G\left(
x,0\right) =\omega ^{\prime }\left( 0\right) $. It follows from diagram (B)
that for all $t\in [0,1]$ 
\begin{equation*}
p_2\circ \omega ^{\prime }(t)=p_2\circ p_1\circ G\left( x,t\right) =p\circ
G(x,t)=p_2\circ F\left( x,t\right) =p_2\circ \omega \left( t\right) \text{.}
\end{equation*}
From this and the condition that $p_2$ has unique path lifting, it follows
that for all $t\in [0,1]$ 
\begin{equation*}
F\left( x,t\right) =\omega \left( t\right) =\omega ^{\prime }\left( t\right)
=p_1\circ G\left( x,t\right) \text{.}
\end{equation*}
Since $x$ was any point of $X$ , we have that $F=p_1\circ G$ and setting $
F^{\prime }=G$ gives us the map we need to complete diagram (A),
demonstrating that $p_1$ is a fibration.\hfill$\square $

\begin{lemma}
Let $X_\infty =\lim\limits_{\leftharpoondown }\{f_i^j,X_i\}$ be the inverse
limit of an inverse sequence for which all the bonding maps $f_i^j$ $(i\leq
j)$ have unique path lifting. Then the projection onto the first coordinate $
f_1:X_\infty \rightarrow X_1;$\ $\left( x_1,x_2,...\right) \mapsto x_1$ has
unique path lifting.
\end{lemma}

\textbf{Proof}: Given paths $\omega $ and $\omega ^{\prime }$ in $X_\infty $
such that $f_1\circ \omega =$ $f_1\circ \omega ^{\prime }$ and $\omega
\left( 0\right) =\omega ^{\prime }\left( 0\right) $, suppose that for some $
t\in [0,1]$ we have $\omega \left( t\right) \neq \omega ^{\prime }\left(
t\right) $. Then for some $n\in \mathbb{N}$, $\omega \left( t\right) $ and $
\omega ^{\prime }\left( t\right) $ disagree on the $n^{th}$ coordinate : $
f_n\circ \omega \left( t\right) \neq f_n\circ \omega ^{\prime }\left(
t\right) $. But we also have that 
\begin{equation*}
f_1^n\circ f_n\circ \omega =f_1\circ \omega =f_1\circ \omega ^{\prime
}=f_1^n\circ f_n\circ \omega ^{\prime }
\end{equation*}
and by hypothesis $f_1^n$ has unique path lifting. This combined with $
f_n\circ \omega \left( 0\right) =f_n\circ \omega ^{\prime }\left( 0\right) $
implies that the paths $f_n\circ \omega $ and $f_n\circ \omega ^{\prime }$
in $X_n$ are equal. This contradicts $f_n\circ \omega \left( t\right) \neq
f_n\circ \omega ^{\prime }\left( t\right) $. Therefore, no such $t$ can
exist and $\omega $ $=\omega ^{\prime }$. \hfill$\square $

\begin{theorem}
Let $X_\infty =\lim\limits_{\leftharpoondown }\{f_i^j,X_i\}$ be the inverse
limit of an inverse sequence for which all the bonding maps $f_i^j$ $(i\leq
j)$ have unique path lifting and suppose that we have a map $p:\tilde{X}
\rightarrow X_\infty $ satisfying the condition that $f_1\circ p$ is a
fibration. Then $p$ is a fibration, and if $\tilde{X}$ is path connected, $
p\left( \tilde{X}\right) $ is a path component of $X_\infty $. If each of
the fibers of $f_1\circ $ $p$ is totally disconnected, then $p$ has unique
path lifting.
\end{theorem}

\textbf{Proof}: Since $f_1$ has unique path lifting by the above lemma, 
\textbf{Lemma \ref{Fib} }applies, implying that $p$ is a fibration. \textbf{
\ }So when $\tilde{X}$ is path connected, $p\left( \tilde{X}\right) $ is a
path component of $X_\infty $ [S; 2.3.1].

Also, if each of the fibers of $f_1\circ p$ is totally disconnected, then
for any $b\in X_\infty $ we have$\;$
\begin{equation*}
p^{-1}(b)\subset p^{-1}\left( f_1^{-1}\left( f_1\left( b\right) \right)
\right) =\left( f_1\circ p\right) ^{-1}\left( f_1\left( b\right) \right)
\end{equation*}
and $\left( f_1\circ p\right) ^{-1}\left( f_1\left( b\right) \right) $ is a
totally disconnected set (see \cite{S}; 2.2.5).\hfill$\square $

Notice that we are not requiring our spaces $X_i$ to be groups and the above
theorem could be extended to include general inverse limits [not just the
inverse limits of inverse sequences].

\begin{corollary}
\label{Pi}For $\kappa \in \mathbb{N}\cup \left\{ \infty \right\} $, let $\pi _{
\overline{M}}:\mathbb{R}^\kappa \rightarrow \sum_{\overline{M}}$ be as in 
\textbf{Definition \ref{Pi def}. }Then $\pi _{\overline{M}}$ is a fibration
with unique path lifting onto $\mathbf{C}_{\overline{M}}$.\hfill$\square $
\end{corollary}

\section{Linear Flows on $\kappa -$Solenoids}

\begin{definition}
For $\mathbf{\omega }=(\omega _1,...,\omega _j,...)$ $\in \mathbb{R}^\kappa -\{
\mathbf{0}\}$, we define $i^{\mathbf{\omega }}:$ $\mathbb{R}$ $\rightarrow \mathbb{
R}^\kappa $ by $t\longmapsto (t\omega _1,...,t\omega _j,...)\overset{def}{=}t
\mathbf{\omega }$.
\end{definition}

Notice that $i^{\mathbf{\omega }}$ topologically embeds $\mathbb{R}$ as a
subgroup of $\mathbb{R}^\kappa $.

\begin{definition}
A \underline{flow} on the space $X$ $\,$is a map $\phi :\mathbb{R\times }
X\rightarrow X$ satisfying the following conditions

\begin{enumerate}
\item  $\phi \left( 0,x\right) =x$ for all $x\in X$

\item  $\phi \left( s,\phi \left( t,x\right) \right) =\phi \left(
s+t,x\right) $ for all $s,t\in \mathbb{R}$.
\end{enumerate}
\end{definition}

\begin{definition}
We define the family of \underline{linear flows} on the $\kappa $-solenoid $
\sum_{\overline{M}}$ $\;$ 
\begin{equation*}
\mathcal{F}_{\overline{M}}=\{\Phi _{\overline{M}}^{\mathbf{\omega }}\mid 
\mathbf{\omega }\in \mathbb{R}^\kappa -\{\mathbf{0}\}\}\text{ to be given by}
\end{equation*}
\begin{equation*}
\Phi _{\overline{M}}^{\mathbf{\omega }}:\mathbb{R\times }\sum\nolimits_{
\overline{M}}\overset{(i^{\mathbf{\omega }},id)}{\longrightarrow }\mathbb{R}
^\kappa \mathbb{\times }\sum\nolimits_{\overline{M}}\overset{(\pi _{\overline{M}
},id)}{\longrightarrow }\mathbf{C}_{\overline{M}}\times \sum\nolimits_{
\overline{M}}\overset{+}{\rightarrow }\sum\nolimits_{\overline{M}}\text{.}
\end{equation*}
\end{definition}

\noindent  It follows directly that each $\Phi _{\overline{M}}^{\mathbf{
\omega }}$ is indeed a flow. Notice that any time$-t$ map $\Phi _{\overline{M
}}^{\mathbf{\omega }}\left( t,\_\right) $ is simply translation by $\pi _{
\overline{M}}\left( t\omega \right) $. This family of flows yields isotopies
between $id_{\sum\nolimits_{\overline{M}}}$ and the translations by elements
of $\mathbf{C}_{\overline{M}}$, and if we replace $\mathbf{C}_{\overline{M}}$
with the path component of $y\in \sum\nolimits_{\overline{M}}-\mathbf{C}_{
\overline{M}}$ and $\pi _{\overline{M}}$ by $y+\pi _{\overline{M}}$ in the
definition of $\Phi _{\overline{M}}^\omega $, we obtain isotopies between
the translations of elements in that path component, but these isotopies are
not flows.

\begin{definition}
The countable set of real numbers $\{\omega _1,...,\omega _i,...\}$ is 
\underline{rationally } \underline{independent} if : 
\begin{equation*}
\left[ k_1\omega _{i_1}^{}+\cdots +k_s\omega _{i_s}^{}=0\text{ for integers }
k_1,....,k_s\text{ (}s\text{ finite)}\right] \Rightarrow \left[ k_1=\cdots
=k_s=0\right] ,
\end{equation*}
and in this case $\mathbf{\omega =}\left( \omega _1,...,\omega _i,...\right) 
$ and the linear flow $\Phi _{\overline{M}}^{\mathbf{\omega }}$ are 
\underline{irrational}.
\end{definition}

\begin{lemma}
\label{Rat}If $\{\omega _1,...,\omega _n\}$ is rationally independent and $N$
is an $n\times n$ invertible matrix with rational entries, then $\{\omega
_1^{\prime },...,\omega _n^{\prime }\}$ is rationally independent, where 
\begin{equation*}
\left( 
\begin{array}{c}
\omega _1^{\prime } \\ 
\vdots \\ 
\omega _n^{\prime }
\end{array}
\right) =N\left( 
\begin{array}{c}
\omega _1 \\ 
\vdots \\ 
\omega _n
\end{array}
\right) \text{.}
\end{equation*}
\end{lemma}

\textbf{Proof}: Suppose that for the integers $k_1,...k_n$ we have $
k_1\omega _1^{\prime }+\cdots +k_n\omega _n^{\prime }=0$. Then 
\begin{equation*}
\left( 
\begin{array}{ccc}
k_1 & \cdots & k_n
\end{array}
\right) \left( 
\begin{array}{c}
\omega _1^{\prime } \\ 
\vdots \\ 
\omega _n^{\prime }
\end{array}
\right) =\left( 
\begin{array}{ccc}
k_1 & \cdots & k_n
\end{array}
\right) N\left( 
\begin{array}{c}
\omega _1 \\ 
\vdots \\ 
\omega _n
\end{array}
\right) =\left( 
\begin{array}{ccc}
q_1 & \cdots & q_n
\end{array}
\right) \left( 
\begin{array}{c}
\omega _1 \\ 
\vdots \\ 
\omega _n
\end{array}
\right) =0
\end{equation*}
where $q_1,\cdots ,q_n$ are rational numbers. The rational independence of $
\{\omega _1,...,\omega _n\}$ then implies that $q_1=\cdots =q_n=0$. Thus, 
\begin{equation*}
\left( 
\begin{array}{ccc}
k_1 & \cdots & k_n
\end{array}
\right) N=\left( 
\begin{array}{ccc}
0 & \cdots & 0
\end{array}
\right) \Rightarrow N^T\left( 
\begin{array}{c}
k_1 \\ 
\vdots \\ 
k_n
\end{array}
\right) =0
\end{equation*}
and so $\left( 
\begin{array}{c}
k_1 \\ 
\vdots \\ 
k_n
\end{array}
\right) \in \ker N^T=\{\left( 
\begin{array}{c}
0 \\ 
\vdots \\ 
0
\end{array}
\right) \}$ and $\{\omega _1^{\prime },...,\omega _n^{\prime }\}$ is
rationally independent by definition.\hfill$\square $

\begin{definition}
$\Lambda _{\overline{M}}^{\mathbf{\omega }}\overset{def}{=}\Phi _{\overline{M
}}^{\mathbf{\omega }}\left( \mathbb{R\times }\{e_{\overline{M}}\}\right) =\{\pi
_{\overline{M}}\left( t\mathbf{\omega }\right) :t\in \mathbb{R}\}\subset 
\mathbf{C}_{\overline{M}}$.
\end{definition}

Notice that $\Lambda _{\overline{M}}^{\mathbf{\omega }}$ is the trajectory
of $e_{\overline{M}}$ for the flow $\Phi _{\overline{M}}^{\mathbf{\omega }}$
and that $\Lambda _{\overline{M}}^{\mathbf{\omega }}$ is a subgroup of $
\sum\nolimits_{\overline{M}}$ since it is the image of $\mathbb{R}$ under the
homomorphism $t\mapsto \pi _{\overline{M}}\left( t\mathbf{\omega }\right) $.

\begin{lemma}
\label{Dense}If $\omega $ is irrational, $\Lambda _{\overline{M}}^{\mathbf{
\omega }}$ is dense in the $\kappa $-solenoid $\sum\nolimits_{\overline{M}}$.
\end{lemma}

\noindent \textbf{Proof}: Let $\mathbf{x=}\left( \mathbf{x}^i\right)
_{i=1}^\infty $ be any point in $\sum\nolimits_{\overline{M}}$ and let $N$
be any neighborhood of $\mathbf{x}$. We need to show that $N$ contains some
point of $\Lambda _{\overline{M}}^{\mathbf{\omega }}$. Since $\mathcal{B}
=\{f_i^{-1}\left( U\right) :U$ is open in $\mathbf{T}^\kappa ,\;i\in \mathbb{N}
\}$ is a basis for the topology of $\sum\nolimits_{\overline{M}}$, there is
a $j\in \mathbb{N}$ and a neighborhood $U$ of $\mathbf{x}^j$ in $\mathbf{T}
^\kappa $ satisfying: $f_j^{-1}\left( U\right) $ $\mathbf{\subset }$ $N$.
Since $U$ is a neighborhood of $\mathbf{x}^j$ in $\mathbf{T}^\kappa $, there
is an $\varepsilon >0$ such that the ball $B$ of radius $\varepsilon $
centered at $\mathbf{x}^j$ in $\mathbf{T}^\kappa $ is contained in $U$. We
have 2 cases: $\kappa =n<\infty $ and $\kappa =\infty $ and we treat the
second case; the first case may be proved using a simplified version of the
same argument.

So we assume $\kappa =\infty $ and seek a point $\pi _{\overline{M}}\left( t
\mathbf{\omega }\right) $ of $\Lambda _{\overline{M}}^{\mathbf{\omega }}$
satisfying $d_\infty \left( \mathbf{x}^j,f_j\left( \pi _{\overline{M}}\left(
t\mathbf{\omega }\right) \right) \right) \\<\varepsilon $. Such a point will
then be contained in $N\cap \Lambda _{\overline{M}}^{\mathbf{\omega }}$
since $f_j^{-1}\left( B\right) \subset $ $f_j^{-1}\left( U\right) $ $\mathbf{
\subset }$ $N$ . First we choose $m$ so that $\sum_{i=m+1}^\infty \frac
1{2^i}<\frac \varepsilon 2$. Then we represent the map $\left(
M_j^{-1}\times id\right) \circ \cdots \circ \left( M_1^{-1}\times id\right) $
as a map $M\times id$, where $M$ is an invertible $k\times k$ matrix for
some integer $k\geq m$. Then with 
\begin{equation*}
\left( 
\begin{array}{c}
\omega _1^{\prime } \\ 
\vdots \\ 
\omega _k^{\prime } \\ 
\vdots
\end{array}
\right) \overset{def}{=}\left( M\times id\right) \left( 
\begin{array}{c}
\omega _1 \\ 
\vdots \\ 
\omega _k \\ 
\vdots
\end{array}
\right) ,
\end{equation*}
the set $\{\omega _1^{\prime },...,\omega _k^{\prime }\}$ and hence $
\{\omega _1^{\prime },...,\omega _m^{\prime }\}$ is rationally independent
by \textbf{Lemma \ref{Rat}}. Kronecker's theorem (see, e.g., \cite{HW}, Thm 444)
then yields integers $p_1,...,p_n$ and a real number $t$ which satisfy the
following system of inequalities: 
\begin{equation*}
\left| t\omega _1^{\prime }-p_1-x_1^j\right| <\frac \varepsilon 2,...,\left|
t\omega _m^{\prime }-p_m-x_m^j\right| <\frac \varepsilon 2\text{,}
\end{equation*}
where for $i=1,...,m$ \ $x_i^j$ (and hence $p_i+x_i^j$) are representatives
in $\mathbb{R}$ for the coordinates of $\mathbf{x}^j$. Then $d_\infty ^{}\left( 
\mathbf{x}^j,f_j\left( \pi _{\overline{M}}\left( t\mathbf{\omega }\right)
\right) \right) <\varepsilon $ as required.\hfill$\square $

\begin{definition}
The flow $\phi :\mathbb{R\times }X\rightarrow X$ is \underline{equivalent} to
the flow $\psi :\mathbb{R\times }Y\rightarrow Y$ if there is a homomorphism $
\alpha :\mathbb{R}\rightarrow \mathbb{R}$ and a homeomorphism $h:X\rightarrow Y$
such that 
\begin{equation*}
\begin{array}{ccc}
\mathbb{R\times }X & \overset{\phi }{\rightarrow } & X \\ 
\downarrow \alpha \times h &  & \downarrow h \\ 
\mathbb{R\times }Y & \overset{\psi }{\rightarrow } & Y
\end{array}
\end{equation*}
commutes and $\alpha $ is increasing $\left[ \text{I,pp.31-2}\right] $, and
we write 
\begin{equation*}
\alpha \times h:\phi \overset{equiv}{\approx }\psi \text{.}
\end{equation*}
\end{definition}

(This is also sometimes referred to as $C^0$ conjugacy).

\begin{lemma}
\label{trans}Let $\Phi _{\overline{M}}^{\mathbf{\omega }}$ be a linear flow
on $\sum_{\overline{M}}$. If $\tau $ is translation by $x\in \sum_{\overline{
M}}$,\ then $id_{\mathbb{R}}\times \tau :\Phi _{\overline{M}}^{\mathbf{\omega }}
\overset{equiv}{\approx }\Phi _{\overline{M}}^{\mathbf{\omega }}$.
\end{lemma}

\textbf{Proof}: We have 
\begin{eqnarray*}
\tau \circ \Phi _{\overline{M}}^{\mathbf{\omega }}\left( t,y\right) &=&\tau
\left( \pi _{\overline{M}}\left( t\mathbf{\omega }\right) +y\right) =\pi _{
\overline{M}}\left( t\mathbf{\omega }\right) +x+y \\
&=&\Phi _{\overline{M}}^{\mathbf{\omega }}\left( t,x+y\right) =\Phi _{
\overline{M}}^{\mathbf{\omega }}\left[ \left( id_{\mathbb{R}}\times \tau
\right) \left( t,y\right) \right] \text{.}
\end{eqnarray*}
\hfill$\square $

Suppose that $\alpha \times g:\Phi _{\overline{M}}^\omega \overset{equiv}{
\approx }\Phi _{\overline{M}}^{\omega ^{\prime }}$ . Then with $\tau $
defined to be translation by $-g\left( e_{\overline{M}}^{}\right) $, we have
by the above lemma: 
\begin{equation*}
id_{\mathbb{R}}\times \tau :\Phi _{\overline{M}}^{\mathbf{\omega }^{\prime }}
\overset{equiv}{\approx }\Phi _{\overline{M}}^{\mathbf{\omega }^{\prime }}
\end{equation*}
and so with $h=\tau \circ g$ we obtain a homeomorphism $h$ which fixes $e_{
\overline{M}}^{}$ and 
\begin{equation*}
\alpha \times h=\left( id_{\mathbb{R}}\times \tau \right) \circ \left( \alpha
\times g\right) :\Phi _{\overline{M}}^\omega \overset{equiv}{\approx }\Phi _{
\overline{M}}^{\omega ^{\prime }}\text{.}
\end{equation*}
Thus, to determine the $\overset{equiv}{\approx }$ classes of the family $
\mathcal{F}_{\overline{M}}$ of linear flows on $\sum_{\overline{M}}$ we need
only consider equivalences which are induced by homeomorphisms of $\sum_{
\overline{M}}$ which fix $e_{\overline{M}}$. We shall see that in fact we
need only consider equivalences induced by automorphisms. In order to
calculate the entropy of automorphisms of $n-$solenoids (in our
terminology), Lind and Ward \cite{LW} show that the general automorphism on an $n$
-solenoid can be represented by a matrix in $GL\left( n,\mathbb{Q}\right) $, as
determined by the dual automorphism on the character group of the solenoid.
We provide here a direct representation of automorphisms on $\kappa $
-solenoids by automorphisms of $\mathbb{R}^\kappa $.

\begin{theorem}
\label{hom}Let $h$ be a homomorphism from the $\kappa $-solenoid $\sum_{
\overline{M}}$ to the $\kappa ^{\prime }$-solenoid $\sum_{\overline{
M^{\prime }}}$. Then there is a unique homomorphism $H$ $:$ $\mathbb{R}^\kappa
\rightarrow \mathbb{R}^{\kappa ^{\prime }}$ making the following diagram
commute 
\begin{equation*}
\begin{array}{ccc}
\mathbb{R}^\kappa & \overset{H}{\rightarrow } & \mathbb{R}^{\kappa ^{\prime }} \\ 
^{\pi _{\overline{M}}}\downarrow &  & \downarrow ^{\pi _{\overline{M^{\prime
}}}} \\ 
\sum_{\overline{M}} & \overset{h}{\rightarrow } & \sum_{\overline{M^{\prime }
}}
\end{array}
\text{.}
\end{equation*}
And the function $f:Hom\left( \ \sum_{\overline{M}},\sum_{\overline{
M^{\prime }}}\right) \rightarrow Hom\left( \mathbb{R}^\kappa \ ,\mathbb{R}^{\kappa
^{\prime }}\right) ;$\ $h\mapsto H$ is one$-$to$-$one. And if $h$ is an
automorphism, $H$ is an automorphism.
\end{theorem}

\textbf{Proof}: Let $h$ be a homomorphism from the $\kappa $-solenoid $\sum_{
\overline{M}}$ to the $\kappa ^{\prime }$-solenoid $\sum_{\overline{
M^{\prime }}}$. Then there is a unique map $H:\left( \mathbb{R}^\kappa ,\mathbf{
0}\right) \rightarrow \left( \mathbb{R}^{\kappa ^{\prime }},\mathbf{0}\right) $
making the following diagrams commute 
\begin{equation*}
\begin{array}{ccc}
&  & \left( \mathbb{R}^{\kappa ^{\prime }},\mathbf{0}\right) \\ 
& ^H\nearrow & \downarrow ^{\pi _{\overline{M^{\prime }}}} \\ 
\left( \mathbb{R}^\kappa ,\mathbf{0}\right) & \underset{h\circ \pi _{\overline{M
}}}{\rightarrow } & \left( \sum_{\overline{M^{\prime }}},e_{\overline{
M^{\prime }}}\right)
\end{array}
\text{ or } 
\begin{array}{ccc}
\left( \mathbb{R}^\kappa ,\mathbf{0}\right) & \overset{H}{\rightarrow } & 
\left( \mathbb{R}^{\kappa ^{\prime }},\mathbf{0}\right) \\ 
^{\pi _{\overline{M}}}\downarrow &  & \downarrow ^{\pi _{\overline{M^{\prime
}}}} \\ 
\left( \sum_{\overline{M}},e_{\overline{M}}\right) & \overset{h}{\rightarrow 
} & \left( \sum_{\overline{M^{\prime }}},e_{\overline{M}}\right)
\end{array}
\end{equation*}
[S; 2.4.2]. Then let $x,y\,$ be any points of $\mathbb{R}^\kappa $. Then 
\begin{eqnarray*}
\pi _{\overline{M^{\prime }}}\left( H\left( x\right) +H\left( y\right)
\right) &=&\pi _{\overline{M^{\prime }}}\circ H\left( x\right) +\pi _{
\overline{M^{\prime }}}\circ H\left( y\right) =h\circ \pi _{\overline{M}
}\left( x\right) +h\circ \pi _{\overline{M}}\left( y\right) \\
&=&h\circ \pi _{\overline{M}}\left( x+y\right) =\pi _{\overline{M^{\prime }}
}\left( H\left( x+y\right) \right) \text{.}
\end{eqnarray*}
From this it follows that $H\left( x\right) +H\left( y\right) -H\left(
x+y\right) \in \ker \pi _{\overline{M^{\prime }}}\subset \mathbb{Z}^{\kappa
^{\prime }}$. Define the following map 
\begin{equation*}
\lambda :\mathbb{R}^\kappa \times \mathbb{R}^\kappa \rightarrow \mathbb{R}^{\kappa
^{\prime }};\left( x,y\right) \mapsto H\left( x\right) +H\left( y\right)
-H\left( x+y\right) \text{.}
\end{equation*}
By the above, $\lambda \left( \mathbb{R}^\kappa \times \mathbb{R}^\kappa \right) $
is a connected subset of the totally disconnected set $\mathbb{Z}^{\kappa
^{\prime }}$ and $\lambda \left( \left( \mathbf{0,0}\right) \right) =\mathbf{
0}$. Thus, $\lambda \left( \mathbb{R}^\kappa \times \mathbb{R}^\kappa \right) =
\mathbf{0}$ and $H\left( x\right) +H\left( y\right) =H\left( x+y\right) $
for all $\left( x,y\right) \in \mathbb{R}^\kappa \times \mathbb{R}^\kappa $ and $H$
is a homomorphism.

Suppose then that $h,h^{\prime }\in Hom\left( \ \sum_{\overline{M}},\sum_{
\overline{M^{\prime }}}\right) $ and that $h\neq h^{\prime }$. Since, for
irrational $\mathbf{\omega }$, $\Lambda _{\overline{M}}^{\mathbf{\omega }
}\subset \mathbf{C}_{\overline{M}}$ is dense we have that $\mathbf{C}_{
\overline{M}}$ is dense. Therefore there is an element $\pi _{\overline{M}
}\left( \mathbf{t}\right) \in \mathbf{C}_{\overline{M}}$ such that $h\left(
\pi _{\overline{M}}\left( \mathbf{t}\right) \right) \neq h^{\prime }\left(
\pi _{\overline{M}}\left( \mathbf{t}\right) \right) $ and hence $\pi _{
\overline{M^{\prime }}}\left( H\left( \mathbf{t}\right) \right) \neq \pi _{
\overline{M^{\prime }}}\left( H^{\prime }\left( \mathbf{t}\right) \right) $,
where $H=$ $f\left( h\right) $ and $H^{\prime }=f\left( h^{\prime }\right) $
. Then we must have $H\left( \mathbf{t}\right) \neq H^{\prime }\left( 
\mathbf{t}\right) $, and so $f$ is one$-$to$-$one.

Suppose then that $h\in Aut\left( \sum_{\overline{M}}\right) $ and $f\left(
h\right) =H:\mathbb{R}^\kappa \rightarrow \mathbb{R}^\kappa $ as above. Then we
have $h^{-1}\in Aut\left( \sum_{\overline{M}}\right) $ and the corresponding
endomorphism $f\left( h^{\prime }\right) =H^{\prime }:\mathbb{R}^\kappa
\rightarrow \mathbb{R}^\kappa $. We then have that $H^{\prime }\circ H$ is the
unique lifting of $id\circ \pi _{\overline{M}}$ and $id_{\mathbb{R}^\kappa }$
also provides such a lifting, so we must have $H^{\prime }\circ H=id_{\mathbb{R}
^\kappa }$. Similarly, $H\circ H^{\prime }=id_{\mathbb{R}^\kappa }$ and $H\in
Aut\left( \mathbb{R}^\kappa \right) $.\hfill$\square $

\begin{lemma}
\label{equi}If $\alpha :\mathbb{R\rightarrow R}$ is multiplication by the
positive number $a$ and if $h\in Isomorphism\left( \sum_{\overline{M}
}\,,\sum_{\overline{M^{\prime }}}\right) $ and if $H:\mathbb{R}^\kappa
\rightarrow \mathbb{R}^\kappa $ satisfies $h\circ \pi _{\overline{M}}=\pi _{
\overline{M^{\prime }}}\circ H$ , then 
\begin{equation*}
\alpha \times h:\Phi _{\overline{M}}^{\mathbf{\omega }}\overset{equiv}{
\approx }\Phi _{\overline{M^{\prime }}}^{\mathbf{\omega }^{\prime }}
\end{equation*}
where $\mathbf{\omega }^{\prime }=\frac 1aH\left( \mathbf{\omega }\right)
\in \mathbb{R}^\kappa $.
\end{lemma}

\textbf{Proof}: With $\mathbf{\omega }^{\prime }=\frac 1aH\left( \mathbf{
\omega }\right) $ we have: 
\begin{eqnarray*}
h\circ \Phi _{\overline{M}}^{\mathbf{\omega }}\left( t,x\right) &=&h\left(
\pi _{\overline{M}}\left( t\mathbf{\omega }\right) +x\right) =h\left( \pi _{
\overline{M}}\left( t\mathbf{\omega }\right) \ \right) +h\left( x\right)
=\pi _{\overline{M^{\prime }}}\left( H\left( t\mathbf{\omega }\right)
\right) +h\left( x\right) \\
&=&\pi _{\overline{M^{\prime }}}\left( at\mathbf{\omega }^{\prime }\right)
+h\left( x\right) =\Phi _{\overline{M^{\prime }}}^{\mathbf{\omega }^{\prime
}}\left( at,h\left( x\right) \right) =\Phi _{\overline{M^{\prime }}}^{
\mathbf{\omega }^{\prime }}\circ \left( \alpha \times h\right) \left(
t,x\right)
\end{eqnarray*}
\hfill$\square $

\begin{lemma}
\label{phase}If the subsets $\Lambda _{\overline{M}}^{\mathbf{\omega }}$ and 
$\Lambda _{\overline{M}}^{\mathbf{\omega }^{\prime }}$ of the $\kappa $
-solenoid $\sum_{\overline{M}}$ are equal, then $\mathbf{\omega =}a\mathbf{
\omega }^{\prime }$ for some $a\in \mathbb{R-\{}0\mathbb{\}}$. And when $\mathbf{
\omega =}a\mathbf{\omega }^{\prime }$, with $S\overset{def}{=}id_{\sum_{
\overline{M}}}$ \thinspace if $a>0\,$ and $S\overset{def}{=}$ $the$ $map$ $
x\mapsto -x$ if $a<0$, we have 
\begin{equation*}
\left( \alpha \times S\right) :\Phi _{\overline{M}}^{\mathbf{\omega }}
\overset{equiv}{\approx }\Phi _{\overline{M}}^{\mathbf{\omega }^{\prime }}
\text{, where }\alpha \text{ is multiplication by }\left| a\right| \text{.}
\end{equation*}
\end{lemma}

\textbf{Proof}: Suppose $\Lambda _{\overline{M}}^{\mathbf{\omega }}=\Lambda
_{\overline{M}}^{\mathbf{\omega }^{\prime }}\overset{def}{=}\Lambda $ and
let $L^{\mathbf{\omega }}=i^\omega \left( \mathbb{R}\right) \subset \mathbb{R}
^\kappa $ and $L^{\mathbf{\omega }^{\prime }}=i^{\omega ^{\prime }}\left( 
\mathbb{R}\right) \subset \mathbb{R}^\kappa $. Since $\Lambda $ is not a single
point and any non-degenerate orbit is a locally one-to-one map from $\mathbb{R}$
, there is an arc $J=\left[ 0,x\right] \subset \mathbb{R}$ so that $p\overset{
def}{=}\pi _{\overline{M}}\circ i^{\mathbf{\omega \ }}\mid _J$ maps $J$
homeomorphically onto its image in $\sum_{\overline{M}}$. Then there is also
an interval $J^{\prime }\subset \mathbb{R}$ having $0$ as an endpoint such that 
$p^{\prime }\overset{def}{=}\pi _{\overline{M}}\circ i^{\mathbf{\omega }
^{\prime }\mathbf{\ }}\mid _{J^{\prime }}$ maps $J^{\prime }$
homeomorphically onto $p\left( J\right) =p^{\prime }\left( J^{\prime
}\right) $, and so there is a homeomorphism $h:\left( J,0\right) \rightarrow
\left( J^{\prime },0\right) $ such that $p^{\prime }\circ h=p$. Then 
\begin{equation*}
\begin{array}{ccc}
\left( J^{\prime },0\right) & \overset{i^{\omega ^{\prime }}}{
\longrightarrow } & \left( \mathbb{R}^\kappa ,\mathbf{0}\right) \\ 
h\uparrow & \searrow ^{p^{\prime }} & \downarrow \pi _{\overline{M}} \\ 
\left( J,0\right) & \underset{p}{\longrightarrow } & \left( \sum_{\overline{M
}},e_{\overline{M}}\right)
\end{array}
,
\end{equation*}
and $i^{\omega ^{\prime }}\circ h$ is the unique lift $\left( J,0\right)
\rightarrow \left( \mathbb{R}^\kappa ,\mathbf{0}\right) $ of $p$. But $i^\omega
\mid _J$ also provides such a lift, and so $i^{\omega ^{\prime }}\circ
h=i^\omega \mid _J$ and $i^\omega \left( J\right) \subset L^{\mathbf{\omega }
^{\prime }}$, from which it follows that $L^{\mathbf{\omega }}=L^{\mathbf{
\omega }^{\prime }}$ and $a\mathbf{\omega }^{\prime }\mathbf{=\omega }$ for
some $a\mathbf{\in }\mathbb{R}$ as claimed.

Suppose then that $\mathbf{\omega =}a\mathbf{\omega }^{\prime }$. Then we
have 
\begin{eqnarray*}
S\circ \Phi _{\overline{M}}^{\mathbf{\omega }}\left( t,x\right) &=&S\left(
\pi _{\overline{M}}\left( t\mathbf{\omega }\right) +x\right) =S\left( \pi _{
\overline{M}}\left( t\mathbf{\omega }\right) \ \right) +S\left( x\right)
=\pi _{\overline{M}}\left( \left| a\right| t\mathbf{\omega }^{\prime
}\right) +S\left( x\right) \\
&=&\Phi _{\overline{M}}^{\mathbf{\omega }^{\prime }}\left( \left| a\right|
t,S\left( x\right) \right) =\Phi _{\overline{M}}^{\mathbf{\omega }^{\prime
}}\circ \left( \alpha \times S\right) \left( t,x\right)
\end{eqnarray*}
\hfill$\square $

Thus, two linear flows $\Phi _{\overline{M}}^{\mathbf{\omega }}$ and $\Phi _{
\overline{M}}^{\mathbf{\omega }^{\prime }}$ whose trajectories determine the
same decomposition of $\sum_{\overline{M}}$ (i.e., linear flows with the
same phase portrait) are equivalent since then the trajectories of $e_{
\overline{M}}$ ($\Lambda _{\overline{M}}^{\mathbf{\omega }}$ and $\Lambda _{
\overline{M}}^{\mathbf{\omega }^{\prime }}$ respectively) are equal.

\begin{definition}
The homeomorphism $h:X\rightarrow Y$ provides a \underline{topological} 
\underline{equivalence} between the flow $\phi :\mathbb{R\times }X\rightarrow X$
and the flow $\psi :\mathbb{R\times }Y\rightarrow Y$ if $h$ maps each
trajectory of $\phi $ onto a trajectory of $\psi $ and if $h$ preserves the
orientation of orbits; that is, for each $x\in X$ there is an increasing
homeomorphism $\alpha _x:\mathbb{R\rightarrow R}$ such that, $h\circ \phi
\left( t,x\right) =\psi \left( \alpha _x\left( t\right) ,h\left( x\right)
\right) $ for all $t\in \mathbb{R}$ $\left[ \text{I, p.32}\right] $. Such $\phi 
$ and $\psi $ are said to be \underline{topologically} \underline{equivalent}
and we write 
\begin{equation*}
h:\phi \overset{top}{\approx }\psi \text{.}
\end{equation*}
\end{definition}

We proceed to determine the $\overset{top}{\approx }$ classes of the
families $\mathcal{F}_{\overline{M}}$, and in the process we shall see that
these $\overset{top}{\approx }$ classes coincide with the $\overset{equiv}{
\approx }$ classes of these families. \textbf{Lemma \ref{trans} }implies
that we need only consider homeomorphisms $h:\left( \sum_{\overline{M}},e_{
\overline{M}}\right) \rightarrow \left( \sum_{\overline{M}},e_{\overline{M}
}\right) $ since topological equivalence is more general than equivalence.
We shall need a result from \cite{Sch}.

\begin{theorem}
Let $G$ be a compact connected topological group, and let $H$ be a locally
compact abelian topological group. Then every $f\in C_e(G,H)=\{maps$ $
G\rightarrow H$ $mapping$ $the$ $identity$ $of$ $G$ $to$ $the\,$ $identity$ $
of$ $H\}$ is homotopic to exactly one $h\in Hom\left( G,H\right) $, and the
homotopy can be chosen to preserve the identity (\cite{Sch}, \textbf{Corollary 2 }
of \textbf{Theorem 2})
\end{theorem}

From this we immediately obtain the following corollary.

\begin{corollary}
Let $G\,$ be a compact connected abelian topological group. Then any
homeomorphism $h:\left( G,e\right) \rightarrow \left( G,e\right) $ is
homotopic to exactly one automorphism.
\end{corollary}

\textbf{Proof}: By our hypotheses on $G$, we obtain $\alpha ,\beta \in
Hom\left( G,G\right) $ with $\alpha $ homotopic to $h$ and $\beta $
homotopic to $h^{-1}$. Then $id_G=h\circ h^{-1}$ is homotopic to $\alpha
\circ \beta \in Hom\left( G,G\right) $. But by theorem, there is only one
element of $Hom\left( G,G\right) $ homotopic to $h\circ h^{-1}$ and $
id_G=h\circ h^{-1}\in Hom\left( G,G\right) $. Therefore, $\alpha \circ \beta
=id_G$ and similarly $\beta \circ \alpha =id_G$. Thus, $\alpha $ is an
automorphism whose uniqueness follows from the uniqueness in the above
theorem.\hfill$\square $

We can actually obtain the automorphism homotopic to $h:\left( G,e\right)
\rightarrow \left( G,e\right) $ as above\textbf{\ }in the following way.
Start with the isomorphism $h^{*}$ of the first \u {C}ech cohomology group
of $G\;$ $\breve{H}^1\left( G\right) $ [$\mathbb{Z}$ coefficients] induced by $
h $. This then yields an automorphism $\iota $ of the dual of $G$\ $\;
\widehat{G}\cong \breve{H}^1\left( G\right) $. The automorphism $\widehat{
\iota }$, the map of $\widehat{\widehat{G}}$ dual to $\iota $, yields, via
the automorphism $G\cong $ $\widehat{\widehat{G}}$ given by Pontryagin
duality, an automorphism $\alpha $ of $G$. This automorphism $\alpha $ is
the automorphism guaranteed by the above.

\begin{theorem}
\label{Class}If $\Phi _{\overline{M}}^{\mathbf{\omega }}\overset{top}{
\approx }\Phi _{\overline{M}}^{\mathbf{\omega }^{\prime }}$, there is an $
a\in Aut\left( \sum_{\overline{M}}\right) $ with $\left( \beta \times
a\right) :\Phi _{\overline{M}}^{\mathbf{\omega }}\overset{equiv}{\approx }
\Phi _{\overline{M}}^{\mathbf{\omega }^{\prime }}$ .
\end{theorem}

\textbf{Proof}: Suppose that $h:\left( \sum_{\overline{M}},e_{\overline{M}
}\right) \rightarrow \left( \sum_{\overline{M}},e_{\overline{M}}\right) $
and $h:\Phi _{\overline{M}}^{\mathbf{\omega }}\overset{top}{\approx }\Phi _{
\overline{M}}^{\mathbf{\omega }^{\prime }}$. Then let $a$ be the unique
automorphism of $\sum_{\overline{M}}$ homotopic to $h$ guaranteed by the
above. Let $F:\left( \sum_{\overline{M}},e_{\overline{M}}\right) \times
[0,1]\rightarrow \left( \sum_{\overline{M}},e_{\overline{M}}\right) $ be a
homotopy of $h$ and $a$ with $F_t\overset{def}{=}F\left( \_,t\right) $. Then
the map $F^{\prime }:\left( \sum_{\overline{M}},e_{\overline{M}}\right)
\times [0,1]\rightarrow \left( \sum_{\overline{M}},e_{\overline{M}}\right)
;\left( x,t\right) \mapsto F_t\left( x\right) -h(x)$ provides a homotopy
between $a-h$ and the constant map $c:\left( \sum_{\overline{M}},e_{
\overline{M}}\right) \rightarrow \{e_{\overline{M}}\}$. Since $\pi _{
\overline{M}}:\mathbb{R}^\kappa \rightarrow \sum_{\overline{M}}$ is a fibration
and $c$ is lifted by the constant map $\sum_{\overline{M}}\rightarrow \{
\mathbf{0}\}$, the map $a-h$ can be lifted by a map $g:\left( \sum_{
\overline{M}},e_{\overline{M}}\right) \rightarrow \left( \mathbb{R}^\kappa ,
\mathbf{0}\right) $ so that the following diagram commutes 
\begin{equation*}
\begin{array}{ccc}
&  & \left( \mathbb{R}^\kappa ,\mathbf{0}\right) \\ 
& ^g\nearrow & \downarrow ^{\pi _{\overline{M}}} \\ 
\left( \sum_{\overline{M}},e_{\overline{M}}\right) & \underset{a-h}{
\longrightarrow } & \left( \sum_{\overline{M}},e_{\overline{M}}\right)
\end{array}
\text{.}
\end{equation*}

Let $A\in Aut\left( \mathbb{R}^\kappa \right) $ be the map which satisfies $\pi
_{\overline{M}}\circ A$ $=a\circ \pi _{\overline{M}}$. Then with $\mathbf{
\varpi =}A\left( \mathbf{\omega }\right) $, we have 
\begin{equation*}
a\left( \pi _{\overline{M}}\left( t\mathbf{\omega }\right) \right) =\pi _{
\overline{M}}\left( A\left( t\mathbf{\omega }\right) \right) =\pi _{
\overline{M}}\left( t\mathbf{\varpi }\right) \text{,}
\end{equation*}
and so $a\left( \Lambda _{\overline{M}}^{\mathbf{\omega }}\right) =\Lambda _{
\overline{M}}^{\mathbf{\varpi }}$. If $\Lambda _{\overline{M}}^{\mathbf{
\varpi }}=\Lambda _{\overline{M}}^{\mathbf{\omega }^{\prime }}$, then 
\textbf{Lemma \ref{phase}} gives us the desired result.

Suppose then that $\Lambda _{\overline{M}}^{\mathbf{\varpi }}\neq \Lambda _{
\overline{M}}^{\mathbf{\omega }^{\prime }}$. Now we have $h:\Phi _{\overline{
M}}^{\mathbf{\omega }}\overset{top}{\approx }\Phi _{\overline{M}}^{\mathbf{
\omega }^{\prime }}$ and so there is an increasing homeomorphism $\alpha :
\mathbb{R\rightarrow R}$ such that, $h\left( \pi _{\overline{M}}\left( t\mathbf{
\omega }\right) \right) =h\circ \Phi _{\overline{M}}^{\mathbf{\omega }
}\left( t,e_{\overline{M}}\right) =\Phi _{\overline{M}}^{\mathbf{\omega }
^{\prime }}\left( \alpha \left( t\right) ,e_{\overline{M}}\right) =\pi _{
\overline{M}}\left( \alpha \left( t\right) \mathbf{\omega }^{\prime }\right) 
$ for all $t\in \mathbb{R}$, and since $h\left( e_{\overline{M}}\right) =e_{
\overline{M}}$ we have $e_{\overline{M}}=\pi _{\overline{M}}\left( \alpha
\left( 0\right) \mathbf{\omega }^{\prime }\right) $. Now for $n\in \mathbb{N}$
define the path 
\begin{equation*}
p_n:[0,n]\overset{i^{\mathbf{\omega }}}{\longrightarrow }\mathbb{R}^\kappa 
\overset{\pi _{\overline{M}}}{\longrightarrow }\sum\nolimits_{\overline{M}}
\overset{\left( a-h\right) }{\longrightarrow }\sum\nolimits_{\overline{M}}
\text{.}
\end{equation*}
Then $\pi _{\overline{M}}\circ g\circ \pi _{\overline{M}}\circ i^{\mathbf{
\omega }}\mid _{[0,n]}=p_n$ and so $g\circ \pi _{\overline{M}}\circ i^{
\mathbf{\omega }}\mid _{[0,n]}$ provides the unique lift $\left(
[0,n],0\right) \rightarrow \left( \mathbb{R}^\kappa ,\mathbf{0}\right) $ of the
path $p_n$. But with 
\begin{equation*}
\gamma :\left( [0,n],0\right) \rightarrow \left( \mathbb{R}^\kappa ,\mathbf{0}
\right) ;\text{\ }t\mapsto t\mathbf{\varpi -}\alpha \left( t\right) \mathbf{
\omega }^{\prime }+\alpha \left( 0\right) \mathbf{\omega }^{\prime }\text{,}
\end{equation*}
we have 
\begin{eqnarray*}
\pi _{\overline{M}}\circ \gamma \left( t\right) &=&\pi _{\overline{M}}\left(
t\mathbf{\varpi -}\alpha \left( t\right) \mathbf{\omega }^{\prime }+\alpha
\left( 0\right) \mathbf{\omega }^{\prime }\right) =\pi _{\overline{M}}\left(
t\mathbf{\varpi }\right) -\pi _{\overline{M}}\left( \alpha \left( t\right) 
\mathbf{\omega }^{\prime }\right) \\
&=&a\left( \pi _{\overline{M}}\left( t\mathbf{\omega }\right) \right)
-h\left( \pi _{\overline{M}}\left( t\mathbf{\omega }\right) \right) =\left(
a-h\right) \circ \left( \pi _{\overline{M}}\circ i^{\mathbf{\omega }}\right)
\left( t\right) =p_n\left( t\right) \text{.}
\end{eqnarray*}
Thus, the uniqueness of $g\circ \pi _{\overline{M}}\circ i^{\mathbf{\omega }
}\mid _{[0,n]}$ yields $\gamma =g\circ \pi _{\overline{M}}\circ i^{\mathbf{
\omega }}\mid _{[0,n]}$ and $g\left( \pi _{\overline{M}}\left( n\mathbf{
\omega }\right) \right) =n\mathbf{\varpi -}\alpha \left( n\right) \mathbf{
\omega }^{\prime }+\alpha \left( 0\right) \mathbf{\omega }^{\prime }$. Let $
\{\pi _{\overline{M}}\left( n_j\mathbf{\omega }\right) \}_{j=1}^\infty
\rightarrow x$ be a convergent subsequence of the sequence $\{\pi _{
\overline{M}}\left( n\mathbf{\omega }\right) \}_{n=1}^\infty $.

With $L^{\mathbf{\varpi }}=\{t\mathbf{\varpi \in }\mathbb{R}^\kappa \mid t
\mathbf{\in }\mathbb{R\}}$ and $L^{\mathbf{\omega }^{\prime }}=\{t\mathbf{
\omega }^{\prime }\mathbf{\in }\mathbb{R}^\kappa \mid t\mathbf{\in }\mathbb{R\}}$,
we have $L^{\mathbf{\varpi }}\neq L^{\mathbf{\omega }^{\prime }}$ [if $L^{
\mathbf{\varpi }}=L^{\mathbf{\omega }^{\prime }}$ there is an $a\mathbf{\in }
\mathbb{R}$ such that $a\mathbf{\omega }^{\prime }\mathbf{=\varpi }$ and $
\Lambda _{\overline{M}}^{\mathbf{\varpi }}=\Lambda _{\overline{M}}^{\mathbf{
\omega }^{\prime }}$]. Then with $L$ $\overset{def}{=}L^{\mathbf{\varpi }
}\oplus L^{\mathbf{\omega }^{\prime }}$, a closed two-dimensional vector
subspace of $\mathbb{R}^\kappa $, we define the topological isomorphism $
\lambda :\mathbb{R}^2\rightarrow L;$\ $\left( s,t\right) \mapsto s\mathbf{
\omega }+t\mathbf{\omega }^{\prime }$ (see, e.g., \cite{Scha}, Chapt1\S 3). Notice
that 
\begin{equation*}
g\left( x\right) \in \overline{\{g\circ \pi _{\overline{M}}\left( n_j\mathbf{
\omega }\right) \}_{j=1}^\infty }=\overline{\{n_j\mathbf{\varpi -}\alpha
\left( n_j\right) \mathbf{\omega }^{\prime }+\alpha \left( 0\right) \mathbf{
\omega }^{\prime }\}_{j=1}^\infty }\subset L.
\end{equation*}
We then have in $\mathbb{R}^2$ 
\begin{eqnarray*}
\lambda ^{-1}\circ g\left( x\right) &=&\lim\limits_j\{\lambda ^{-1}\circ
g\left( \pi _{\overline{M}}\left( n_j\mathbf{\omega }\right) \right) \} \\
&=&\lim\limits_j\{\lambda ^{-1}\left( n_j\mathbf{\varpi -}\alpha \left(
n_j\right) \mathbf{\omega }^{\prime }+\alpha \left( 0\right) \mathbf{\omega }
^{\prime }\right) \}=\lim\limits_j\{\left( n_j,\mathbf{-}\alpha \left(
n_j\right) +\alpha \left( 0\right) \right) \},
\end{eqnarray*}
which is impossible since $\{\left( n_j,\mathbf{-}\alpha \left( n_j\right)
+\alpha \left( 0\right) \right) \}$ is unbounded. We must therefore have $
\Lambda _{\overline{M}}^{\mathbf{\varpi }}=\Lambda _{\overline{M}}^{\mathbf{
\omega }^{\prime }}$.\hfill$\square $

This reduces the problem of determining the $\overset{equiv}{\approx }$ and $
\overset{top}{\approx }$ classes of the families $\mathcal{F}_{\overline{M}}$
to determining the image of \ $f:Aut\left( \sum_{\overline{M}}\right)
\rightarrow Aut\left( \mathbb{R}^\kappa \right) $ as in \textbf{Theorem \ref
{hom}}. Two linear flows $\Phi _{\overline{M}}^\omega $ and $\Phi _{
\overline{M}}^{\omega ^{\prime }}$ are (topologically) equivalent if and
only if there is an $h\in Aut\left( \sum_{\overline{M}}\right) $ whose lift $
H$ to an automorphism of $\mathbb{R}^\kappa $ satisfies: $a\omega ^{\prime
}=H\left( \omega \right) $ for some $a\in \mathbb{R-}\left\{ 0\right\} $.
Generally, if $\Phi _{\overline{M}}^\omega $ $\overset{equiv}{\approx }$ $
\Phi _{\overline{M}}^{\omega ^{\prime }}$ and the rank of the subgroup of $
\left( \mathbb{R},\mathbb{+}\right) $ generated by $\left\{ \omega _1,\omega
_2,...\right\} $ is $\rho $, then the rank of the subgroup generated by $
\left\{ \omega _1^{\prime },\omega _2^{\prime },...\right\} $ is also $\rho $
: the closure of the corresponding trajectory in each case is a $\rho -$
solenoid. However, it is important to realize that the image $f$ $\left(
Aut\left( \sum_{\overline{M}}\right) \right) $ depends on the $\kappa $
-solenoid $\sum_{\overline{M}}$ and so which flows of the same rank are
equivalent depends on $\sum_{\overline{M}}$. Next we shall give a specific
classification of these automorphisms on $n$-solenoids whose bonding maps
are all represented by diagonal matrices, which correspond to the finite
product of $1$-solenoids.

\section{Classifying Automorphisms on the Finite Product of 1-Solenoids}

\subsection{Comparing $1-$Solenoids}

For a sequence of non-zero integers $P=(p_1,p_2,...)$, we have the
corresponding $1$-solenoid $\sum\nolimits_P$ where the bonding map $
f_i^{i+1} $ is multiplication by $p_i$ in $\mathbf{T}^1$. [This is
consistent with our established terminology; usually the $p_i$ are required
to be positive primes, but we include the case $p_i=1$].

\begin{definition}
Sequences of non-zero integers $P=(p_1,p_2,...)$ and $Q=(q_1,q_2,...)$ are 
\underline{equivalent}, denoted $P\cong Q$, if and only if $\sum_P$ is
topologically isomorphic with $\sum_Q$. We let $\bar{P}$ denote the $\cong $
class of $P$.
\end{definition}

\begin{definition}
For a given sequence of non-zero integers $P=(p_1,p_2,...)$, $\left|
P\right| \overset{def}{=}(\left| p_1\right| ,\left| p_2\right| ,...)$.
\end{definition}

\begin{proposition}
$P\cong \left| P\right| $.
\end{proposition}

\textbf{Proof}: With $sgn\left( i\right) $ is defined as the map $
S^1\rightarrow S^1$ given by multiplication by $sign\left( p_1\cdots
p_i\right) $, we have the following topological isomorphism $
\sum_P\rightarrow \sum_{\left| P\right| }$ represented by the vertical maps
in the following commutative diagram 
\begin{equation*}
\begin{array}{cccccc}
S^1 & \overset{p_1}{\longleftarrow } & S^1 & \overset{p_2}{\longleftarrow }
& S^1 & \longleftarrow \cdots \\ 
id\downarrow &  & sgn\left( 1\right) \downarrow &  & sgn\left( 2\right)
\downarrow &  \\ 
S^1 & \overset{\left| p_1\right| }{\longleftarrow } & S^1 & \overset{\left|
p_2\right| }{\longleftarrow } & S^1 & \longleftarrow \cdots
\end{array}
\end{equation*}
\hfill$\square $

\begin{definition}
Given a sequence $P=(p_1,p_2,...)$ of non-zero integers, we define the 
\underline{derived sequence} $P^{\prime }=(p_1^{\prime },p_2^{\prime },...)$
to be the equivalent sequence of primes and $1$'s obtained by the prime
factorization of the $\left| p_i\right| $ in sequence, with the factors of $
\left| p_i\right| $ ordered by magnitude $\left[ 1^{\prime }\text{s are left
unchanged}\right] $.
\end{definition}

For example: for $P=(6,1,-90,...),$ $P^{\prime }=(2,3,1,2,3,3,5,...)$.

\begin{definition}
We define the partial order $\leq $ on sequences of non-zero integers: $
P\leq Q$ iff a finite number of terms may be deleted from $P^{\prime }$ so
that each prime occurring in this deleted sequence$\,$occurs in $Q^{\prime }$
with the same or greater cardinality.
\end{definition}

Since $\left( 1,1,...\right) $ is a minimal element of the partial order
which is equivalent with any other minimal element, we can summarize the
classification of all $1$-solenoids in this new terminology as follows. This
classification was conjectured by Bing \cite{B}, while the first proof in print
appears in \cite{McC}$;$ see also \cite{AF}.

\begin{theorem}
$\left( \left[ P\leq Q\right] \text{ and }\left[ Q\leq P\right] \right) $ $
\Leftrightarrow \left[ P\cong Q\right] $.
\end{theorem}

And so $\leq $ induces a partial ordering $\preceq $ on $\cong $ classes.

\begin{definition}
We define the sequence of primes 
\begin{equation*}
\prod =\left( 2,3,2,5,3,2,...\right) =\left(
p_1^1,p_1^2,p_2^1,p_1^3,p_2^2,p_3^1,...\right) =\left( \phi _1,\phi
_2,...\right) ,
\end{equation*}
where $p_i^k$ is the $k^{th}$ prime for all $i$. That is, with\ $f:\mathbb{
N\times N\rightarrow N}$ denoting the standard ``diagonal'' bijection given
by $\left( k,i\right) \mapsto \left. 
\begin{array}{c}
\underline{\left( k+i-2\right) \left( k+i-1\right) } \\ 
2
\end{array}
\right. +i$, $\prod $ is the sequence of primes whose $j^{th}$ member $\phi
_j$ satisfies the following: $\phi _j=p_i^k$ , where $\left( k,i\right) =$ $
f^{-1}\left( j\right) $; equivalently, $\phi _{f\left( k,i\right) }=p_i^k$.
\end{definition}

Notice that under the partial ordering $\preceq $ $\overline{\left(
1,1,...\right) }$ is the minimum element and $\overline{\Pi }$ is the
maximum element.

\begin{definition}
For the $\,$prime $p$ and a sequence of primes and $1$'s $P$, $card_p\left(
P\right) $ is the cardinality with which $p$ occurs in $P$.
\end{definition}

We now introduce a way of arranging finitely many sequences of integers
which is convenient for our purposes.

\begin{definition}
Given sequences of non-zero integers 
\begin{equation*}
P_1=\left( p_1^1,p_2^1,...\right) ,...,P_n=\left( p_1^n,p_2^n,...\right) ,
\end{equation*}
we define their \underline{proper} \underline{arrangement} to be the
sequences 
\begin{equation*}
Q_1=\left( q_1^1,q_2^1,...\right) ,...,Q_n=\left( q_1^n,q_2^n,...\right)
\end{equation*}
defined as follows: let $\kappa \overset{def}{=}\max \{card_2\left(
P_1^{\prime }\right) ,...,card_2\left( P_n^{\prime }\right) \}$ and $m
\overset{def}{=}\min \{i:card_2\left( P_i^{\prime }\right) =\kappa \}$. For
all $i\in \mathbb{N}$ 
\begin{equation*}
q_{f(1,i)}^m\overset{def}{=}\left\{ 
\begin{array}{c}
2\text{ for }i\leq \kappa \\ 
1\text{ otherwise}
\end{array}
\right. ,\left( 
\begin{array}{c}
q_{f(1,i)}^m=2\text{ for all }i\text{ if }\kappa =\infty \\ 
q_{f(1,i)}^m=1\text{ for all }i\text{ if }\kappa =0
\end{array}
\right) .
\end{equation*}
For any $P_j$ satisfying $P_m\leq P_j$ define $q_{f(1,i)}^j=q_{f(1,i)}^m$
for all $i\in \mathbb{N}$. Let $\{$ $j_1,...,j_k\}\subset \{1,...,n\}$ be the
indices $\ell $ for which the terms $q_{f(1,i)}^\ell $ have not been defined
(possibly empty). Repeat the above procedure on $P_{j_1},...,P_{j_k}.$
Repeat the procedure as many times as needed until $q_{f(1,i)}^j$ is defined
for all $j\in \{1,...,n\}$ and all $i\in \mathbb{N}$. Proceed recursively, the
occurrence of the $k^{th}$ prime in the sequences determines the values of $
q_{f(k,i)}^j$ for $\left( j,i\right) \in \{1,...,n\}\times \mathbb{N}$.
\end{definition}

Notice that $Q_i\cong P_i$ since to compose $Q_i$ we add at most a finite
number of prime factors to those of $P_i^{\prime }$.

\begin{definition}
$P_1=\left( p_1^1,p_2^1,...\right) ,...,P_n=\left( p_1^n,p_2^n,...\right) $
are \underline{properly arranged} iff the proper arrangement of $P_1,...,P_n$
is $P_1,...,P_n$.
\end{definition}

\begin{definition}
If $P$ and $Q$ are properly arranged and $P$ $\geq $ $Q$ or if $P=Q$, for $
n\in \mathbb{Z}$ we define the map $n_{P\rightarrow Q}:\sum_P\rightarrow \sum_Q$
to be the map represented by the vertical maps in the following commutative
diagram 
\begin{equation*}
\begin{array}{cccccccc}
S^1 & \overset{p_1}{\longleftarrow } & S^1 & \overset{p_2}{\longleftarrow }
& S^1 & \overset{p_3}{\longleftarrow } & \cdots &  \\ 
n\downarrow &  & nr_1\downarrow &  & nr_2\downarrow &  &  &  \\ 
S^1 & \overset{q_1}{\longleftarrow } & S^1 & \overset{q_2}{\longleftarrow }
& S^1 & \overset{q_3}{\longleftarrow } & \cdots & 
\end{array}
,\text{ }
\end{equation*}
where $r_i= 
\begin{array}{c}
\underline{p_1\cdots p_i} \\ 
q_1\cdots q_i
\end{array}
\in \mathbb{Z-\{}0\}$ since $P$ $\geq $ $Q$. We also denote $n_{P\rightarrow P}$
by $n_P$.
\end{definition}

Then for $n\neq 0$ $n_{P\rightarrow Q}$ is a topological epimorphism since
the maps $n,nr_1,...:S^1\rightarrow S^1$ are topological epimorphisms.

\subsection{Solenoidal Arithmetic}

\begin{definition}
For the prime $p$, we say $p\mid P$ ($p$ divides $P$) if and only if $p$
divides infinitely many $p_i\in P$.
\end{definition}

The following lemma is a translation of the result [K, 2.4] into our
terminology; it may also be shown number theoretically.

\begin{lemma}
\label{prime}For the prime $p\in \mathbb{N}$,\ $p_P$ is a topological
isomorphism $\Leftrightarrow $ $p\mid P$.$\square $
\end{lemma}

\begin{definition}
For $r\in \mathbb{R}$, define $r^{P\rightarrow Q}:\mathbf{C}_P\rightarrow 
\mathbf{C}_Q$ by $r^{P\rightarrow Q}\left( \pi _P(s)\right) =\pi _Q(rs)$ if $
P$ is not eventually $1.$ And if $P$ is eventually $1$ and the product of
all the terms ($\neq 1$) of $P$ is $p$, then $\pi _P:[0,p)\rightarrow 
\mathbf{C}_P$ is one$-$to$-$one and onto, and we define $r^{P\rightarrow Q}:
\mathbf{C}_P\rightarrow \mathbf{C}_Q$ by $r^{P\rightarrow Q}\left( \pi
_P(s)\right) =\pi _Q(rs)$ , $s\in [0,p)$.
\end{definition}

Notice that if $P$ is eventually $1$ and $Q$ is not eventually $1$, $
r^{P\rightarrow Q}$ is not continuous; for if it were, then $
\lim\limits_{s\rightarrow p^{-}}r^{P\rightarrow Q}\left( \pi _P(s)\right) $
would be $\pi _Q(rp)\neq e_Q$, while $\lim\limits_{s\rightarrow
0+}r^{P\rightarrow Q}\left( \pi _P(s)\right) =e_Q$ and $\lim\limits_{s
\rightarrow p^{-}}\pi _P(s)=\lim\limits_{s\rightarrow 0+}\pi _P(s)=e_P$.

\begin{lemma}
If $r\in \mathbb{R-Q}$, then $r^{P\rightarrow Q}$ is not continuous.
\end{lemma}

\textbf{Proof}: If $P$ is eventually $1,$ then $r^{P\rightarrow Q}$ followed
by projection onto the first $S^1$ factor is multiplication by an irrational
on $S^1$, which is not continuous. Suppose then that $P$ is not eventually $
1 $. Now $e_Q^{}\in J\overset{def}{=}$ $\cup _{i\in \mathbb{Z}}\{\pi _Q\left(
(-\frac 14,\frac 14)+iq_1\right) \}=\cup _{i\in \mathbb{Z}}J_i$ is a basic open
subset of $\mathbf{C}_Q$, and if $r^{P\rightarrow Q}$ were continuous, there
would be a basic open subset $I=$ $\cup _{i\in \mathbb{Z}}\{\pi _P\left(
(-\delta ,\delta )+ip_1\cdots p_k\right) \}$ containing $e_P$ such that $
r^{P\rightarrow Q}\left( I\right) =\cup _{i\in \mathbb{Z}}\{\pi _Q\left(
(-\alpha ,\alpha )+irp_1\cdots p_k\right) \}=\cup _{i\in \mathbb{Z}}I_i\subset
J $, where $\alpha =\left| r\delta \right| $. But the centers of successive $
I_i$ have preimages under $\pi _Q$ which are not spaced by an integral
amount as the centers of the $J_i$ are $\Rightarrow $ there is some $I_k$
which is not contained in any $J_i$, contradicting $r^{P\rightarrow Q}\left(
I\right) \subset J$.\hfill $\square $

\begin{lemma}
If $r=\frac cd\in \mathbb{Q-\{}0\}$ and $d$ has a prime factor $p$ which does
not divide $P$ and $\gcd (c,d)=1$, then $r^{P\rightarrow P}$ is not
continuous.
\end{lemma}

\textbf{Proof}: Let $\{s_n\}=\{\pi _P\left( p_1\cdots p_n\right) \}$. Then $
\{s_n\}\rightarrow e_P$. Since $p$ does not divide $P,$ there is an $N$ such
that for all $n\geq N,$ $p$ does not divide $p_n$. For $n\geq N$ the $\left(
N+1\right) ^{th}\,$ coordinate of $r^{P\rightarrow P}\left( s_n\right) = 
\begin{array}{c}
\underline{cp_N\cdots p_n} \\ 
d
\end{array}
=k+\frac \ell d$, where $k\in \mathbb{Z}$ and $\ell \in
\{1,...,d-1\}\Rightarrow $ $d_1^\infty \left( r^{P\rightarrow P}\left(
s_n\right) ,e_P\right) \geq \frac 1{2^{N+1}d}\Rightarrow \{r^{P\rightarrow
P}\left( s_n\right) \}$ does not converge to $e_P$.\hfill $\square $

Notice that if $r=\frac 1d$ and all prime factors of $d$ divide $P$, then $
r^{P\rightarrow P}$ is a topological isomorphism by \textbf{Lemma \ref{prime}
}.

\begin{lemma}
If $r=\frac cd\in \mathbb{Q}-\{0\}$ and $P>Q$ and $d$ has a prime factor $p$
which does not divide $P$, then $r^{P\rightarrow Q}$ is not continuous.
\end{lemma}

\textbf{Proof}: Very similar to the above proof.\hfill $\square $

\begin{lemma}
\label{greater}If $r=\frac cd\in \mathbb{Q}$ and $P>Q$ are properly arranged
and $d$ is the product of primes which divide $P$, then $r^{P\rightarrow Q}$
is continuous.
\end{lemma}

\textbf{Proof}: $r^{P\rightarrow Q}=c_{P\rightarrow Q}\circ \left( \frac
1d\right) ^{P\rightarrow P}$\hfill$\square $

\begin{lemma}
If $r=\frac cd\in \mathbb{Q}-\{0\}$ and if $P<Q$ or $P$ and $Q$ are not
comparable, then $r^{P\rightarrow Q}$ is not continuous.
\end{lemma}

\textbf{Proof}: Let $\{s_n\}=\{\pi _P\left( p_1\cdots p_n\right) \}$. Then $
\{s_n\}\rightarrow e_P$. By our hypothesis, there is an $N$ such that $
q_1\cdots q_N$ does not divide $cp_1\cdots p_m$ for all $m$. Then the $
\left( N+1\right) ^{th}\,$ coordinate of $r^{P\rightarrow Q}\left(
s_n\right) = 
\begin{array}{c}
\underline{cp_1\cdots p_n} \\ 
dq_1\cdots q_N
\end{array}
=k+\frac \ell {dq_1\cdots q_N}$, where $k\in \mathbb{Z}$ and $\ell \in
\{1,...,dq_1\cdots q_N-1\}\Rightarrow $ $\{r^{P\rightarrow Q}\left(
s_n\right) \}$ does not converge to $e_Q$. \hfill$\square $

\begin{definition}
We define $r\in \mathbb{R}$ to be a \underline{proper $P\rightarrow Q$
multiplier} $[$\/or \underline{proper} \underline{multiplier} if the context
is clear\/$]$ if the corresponding function $r^{P\rightarrow Q}$ is a
topological homomorphism.
\end{definition}

Thus, we may summarize our above results as follows:

If $P\geq Q$, then $r$ is a non-zero proper $P\rightarrow Q$ multiplier if
and only if $r=\frac cd\in \mathbb{Q-\{}0\}$ and all prime factors of $d$
divide $P;$ if $P<Q$ or $P$ and $Q$ are not comparable, then there is no
non-zero proper $P\rightarrow Q$ multiplier.

\begin{definition}
When $r$ is a proper $P\rightarrow Q$ multiplier we define the corresponding
map $r_{P\rightarrow Q}:\sum_P\rightarrow \sum_Q$ to be the map which
extends $r^{P\rightarrow Q}$.
\end{definition}

Notice that $r_{P\rightarrow Q}$ is well-defined since we can write it as
the composition of topological isomorphisms $\sum_P\rightarrow \sum_P$
(corresponding to the composition of maps given by the factors in the
denominator) and the epimorphism $n_{P\rightarrow Q}$ [when $n\neq 0$] for
some $n\in \mathbb{Z}$. Thus, we obtain the following result.

\begin{lemma}
$r_{P\rightarrow Q}$ is a topological epimorphism for any non-zero proper
multiplier $r$; if $P=Q$, then $r_{P\rightarrow P}$ is a topological
isomorphism if and only if $r$ and $\frac 1r$ are proper multipliers.\hfill$
\square $
\end{lemma}

\begin{definition}
We define $r\in \mathbb{R}$ to be a \underline{$\left( P\right) \;$
iso-multiplier} if and only if $r$ and $\frac 1r$ are proper $P\rightarrow P$
multipliers.
\end{definition}

\subsection{Cartesian Products of $1-$Solenoids as $n-$Solenoids and their
Automorphisms}

\begin{definition}
For $n$ sequences of non-zero integers $P_1=\left( p_1^1,p_2^1,...\right)
,...,P_n=\left( p_1^n,p_2^n,...\right) $ we define the $n$-solenoid $
\sum\nolimits_{\overline{P}}$ to be the $n$-solenoid corresponding to the
sequence of matrices 
\begin{equation*}
\overline{P}=\left( \left( 
\begin{array}{ccc}
p_1^1 &  & 0 \\ 
& \ddots &  \\ 
0 &  & p_1^n
\end{array}
\right) ,\left( 
\begin{array}{ccc}
p_2^1 &  & 0 \\ 
& \ddots &  \\ 
0 &  & p_2^n
\end{array}
\right) ,...\right) \text{.}
\end{equation*}
\end{definition}

\begin{definition}
$d_{\left( P_1,...,P_n\right) }$ is the metric on the product $
\prod_{i=1}^n\sum\nolimits_{P_i}$ given by 
\begin{equation*}
d_{\left( P_1,...,P_n\right) }\left( \left( \langle x_1^j\rangle
_{j=1}^\infty ,...,\langle x_n^j\rangle _{j=1}^\infty \right) ,\left(
\langle y_1^j\rangle _{j=1}^\infty ,...,\langle y_n^j\rangle _{j=1}^\infty
\right) \right) \overset{def}{=}\,\sum_{i=1}^n\frac 1{2^i}d_\infty ^{}\left(
\langle x_i^j\rangle _{j=1}^\infty ,\langle y_i^j\rangle _{j=1}^\infty
\right) \text{.}
\end{equation*}
\end{definition}

\begin{lemma}
\label{prod}$\prod_{i=1}^n\sum\nolimits_{P_i}$ is topologically isomorphic
to $\sum\nolimits_{\overline{P}}$.
\end{lemma}

\textbf{Proof}: Define $\mathfrak{i}:\prod_{i=1}^n\sum\nolimits_{P_i}\rightarrow
\sum\nolimits_{\overline{P}}$ by 
\begin{equation*}
\left( \langle x_1^j\rangle _{j=1}^\infty ,...,\langle x_n^j\rangle
_{j=1}^\infty \right) \mapsto \left( \langle x_1^j,...,x_n^j\rangle
_{j=1}^\infty \right) \in \sum\nolimits_{\overline{P}}\subset
\prod_{j=1}^\infty \mathbf{T}^n.
\end{equation*}
The function $\mathfrak{i}$ is well-defined since $x_i^j=p_j^ix_i^{j+1}$. By
construction, $\mathfrak{i}$ is an algebraic isomorphism. And $\mathfrak{i}$ is a
homeomorphism; in fact, $\mathfrak{i}$ is an isometry: 
\begin{equation*}
d_{\left( P_1,...,P_n\right) }\left( \left( \langle x_1^j\rangle
_{j=1}^\infty ,...,\langle x_n^j\rangle _{j=1}^\infty \right) ,\left(
\langle y_1^j\rangle _{j=1}^\infty ,...,\langle y_n^j\rangle _{j=1}^\infty
\right) \right) =\sum_{i=1}^n\frac 1{2^i}d_\infty ^{}\left( \langle
x_i^j\rangle _{j=1}^\infty ,\langle y_i^j\rangle _{j=1}^\infty \right)
\end{equation*}
\begin{equation*}
=\sum_{i=1}^n\frac 1{2^i}\left( \sum_{j=1}^\infty \frac 1{2^j}d_1\left(
\left( x_i^j,y_i^j\right) \right) \right) =\sum_{i=1}^n\sum_{j=1}^\infty
\frac 1{2^i2^j}d_1\left( x_i^j,y_i^j\right) \text{ and}
\end{equation*}
\ 
\begin{equation*}
d_n^\infty \left( \mathfrak{i}\left( \langle x_1^j\rangle _{j=1}^\infty
,...,\langle x_n^j\rangle _{j=1}^\infty \right) ,\mathfrak{i}\left( \langle
y_1^j\rangle _{j=1}^\infty ,...,\langle y_n^j\rangle _{j=1}^\infty \right)
\right) =d_n^\infty \left( \langle \mathbf{x}^j\rangle _{j=1}^\infty
,\langle \mathbf{y}^j\rangle _{j=1}^\infty \right) =
\end{equation*}
\begin{equation*}
\sum_{j=1}^\infty \frac 1{2^j}d_n\left( \langle x_1^j,...,x_n^j\rangle
,\langle y_1^j,...,y_n^j\rangle \right) =\sum_{j=1}^\infty \frac
1{2^j}\left( \sum_{i=1}^n\frac 1{2^i}d_1\left( \left( x_i^j,y_i^j\right)
\right) \right)
\end{equation*}
\begin{equation*}
=\sum_{j=1}^\infty \sum_{i=1}^n\frac 1{2^j2^i}d_1\left( x_i^j,y_i^j\right) 
\text{.}
\end{equation*}
\hfill$\square $

Notice that on $\prod_{i=1}^n\mathbf{C}_{P_i}$ $\mathfrak{i}$ takes on the
simple form 
\begin{equation*}
\left( \pi _{P_1}\left( t_1\right) ,...,\pi _{P_n}\left( t_n\right) \right)
\mapsto \pi _{\overline{P}}\left( \left( t_1,...,t_n\right) \right) .
\end{equation*}

\begin{theorem}
\label{aut}If $P_1=\left( p_1^1,p_2^1,...\right) ,...,P_n=\left(
p_1^n,p_2^n,...\right) $ are properly arranged and if $h$ is an automorphism
of $\sum\nolimits_{\overline{P}}$ with the corresponding automorphism $H:
\mathbb{R}^n\rightarrow \mathbb{R}^n$ satisfying $h\circ \pi _{\overline{P}}=\pi _{
\overline{P}}\circ H$ and if $A=\left( a_{ij}\right) $ is the matrix
representing $H$ and if $A^{-1}=\left( b_{ij}\right) $, then the entries $
a_{ij}$ and $b_{ij}$ are proper $P_j\rightarrow P_i$ multipliers for each $
i,j\in \{1,...,n\}$. And if $A=\left( a_{ij}\right) $ is an invertible
matrix with inverse $A^{-1}=\left( b_{ij}\right) $ and if all the entries $
a_{ij}$ and $b_{ij}$ are proper $P_j\rightarrow P_i$ multipliers, then there
is an automorphism of $\sum\nolimits_{\overline{P}}$ represented by $A$. In
the notation of \textbf{Theorem \ref{hom}}, 
\begin{equation*}
f\left( Aut\left( \sum\nolimits_{\overline{P}}\right) \right) =\{A\in
GL\left( n,\mathbb{Q}\right) :\text{the entries of }A\text{ and }A^{-1}\text{
are proper }P_j\rightarrow P_i\text{ multipliers}\}\text{.}
\end{equation*}
\end{theorem}

\textbf{Proof}: Let $h,H$ and $A\,$ be as in the statement of the theorem.
Fix an entry $a_{ij}$ of $A$. With $\phi _k:\prod_{\ell
=1}^n\sum\nolimits_{P_\ell }\rightarrow \sum\nolimits_{P_k}$ projection onto
the $k^{th}$ factor, we have that $\phi _j\circ \mathfrak{i}^{-1}$ maps $R_j
\overset{def}{=}\{\pi _{\overline{P}}\left( \mathbf{t}\right) \in
\sum\nolimits_{\overline{P}}\mid \mathbf{t}=\left( t_1,...,t_n\right) \in 
\mathbb{R}^n$ and $t_\ell =0$ for $\ell \neq j\}$ isomorphically onto $\mathbf{C
}_{P_j}:$ 
\begin{equation*}
\pi _{\overline{P}}\left( \left( 0,...,t,0,..,0\right) \right) \overset{
\mathfrak{i}^{-1}}{\mapsto }\left( \pi _{P_1}\left( 0\right) ,...,\pi
_{P_j}\left( t\right) ,\pi _{P_{j+1}}\left( 0\right) ,...,\pi _{P_n}\left(
0\right) \right) \overset{\phi _j}{\mapsto }\left( \pi _{P_j}\left( t\right)
\right) \text{.}
\end{equation*}
Thus, we have the map $\mu \overset{def}{=}$ $\phi _i\circ \mathfrak{i}
^{-1}\circ h\circ \left( \phi _j\circ \mathfrak{i}^{-1}\right) ^{-1}\mid _{
\mathbf{C}_{P_j}}:\mathbf{C}_{P_j}\rightarrow \mathbf{C}_{P_i}\subset
\sum\nolimits_{P_i}$. And for $t\in \mathbb{R}:$ 
\begin{eqnarray*}
\mu \left( \pi _{P_j}\left( t\right) \right) &=&\phi _i\circ \mathfrak{i}
^{-1}\circ h\left( \pi _{\overline{P}}\left( \left( 0,...,t,0,..,0\right)
\right) \right) =\phi _i\circ \mathfrak{i}^{-1}\circ \pi _{\overline{P}}\circ
H\left( \left( 0,...,t,0,..,0\right) \right) \\
&=&\phi _i\circ \mathfrak{i}^{-1}\circ \pi _{\overline{P}}\left(
a_{1j}t,...,a_{nj}t\right) =\pi _{P_i}\left( a_{ij}t\right) \text{,}
\end{eqnarray*}
and so $a_{ij}$ is a proper $P_j\rightarrow P_i$ multiplier since the map $
\mu $ equals $\left( a_{ij}\right) ^{P_j\rightarrow P_i}$. Similarly, each $
b_{ij}\,$ is a proper $P_j\rightarrow P_i$ multiplier since $A^{-1}$
represents $h^{-1}$.

Given an invertible matrix $A=\left( a_{ij}\right) $ with inverse $
A^{-1}=\left( b_{ij}\right) $where the entries $a_{ij}$ and $b_{ij}$ are
proper $P_j\rightarrow P_i$ multipliers for all $i,j\in \{1,...,n\}$, we
define the map $A_{\overline{P}}$ to be the map on $\sum\nolimits_{\overline{
P}}$ conjugate via $\mathfrak{i}^{-1}$ to the map $\mathcal{A}$ $:$ $\prod_{\ell
=1}^n\sum\nolimits_{P_\ell }\rightarrow \prod_{\ell
=1}^n\sum\nolimits_{P_\ell };$ 
\begin{equation*}
\left( 
\begin{array}{c}
x_1 \\ 
\vdots \\ 
x_n
\end{array}
\right) \overset{\mathcal{A}}{\mapsto }\left( 
\begin{array}{ccc}
\left( a_{11}\right) _{P_1\rightarrow P_1}\left( x_1\right) & +\cdots + & 
\left( a_{11}\right) _{P_n\rightarrow P_1}\left( x_n\right) \\ 
& \vdots &  \\ 
\left( a_{n1}\right) _{P_1\rightarrow P_n}\left( x_1\right) & +\cdots + & 
\left( a_{nn}\right) _{P_n\rightarrow P_n}\left( x_n\right)
\end{array}
\right) \text{; }A_{\overline{P}}\overset{def}{=}\mathfrak{i}\circ \mathcal{A}
\circ \mathfrak{i}^{-1}\text{.}
\end{equation*}
Then $A_{\overline{P}}$ $\circ \pi _{\overline{P}}=\pi _{\overline{P}}\circ
A $ as desired and $\left( A_{\overline{P}}\right) ^{-1}=\left(
A^{-1}\right) _{\overline{P}}.$\hfill$\square $

Any finite product of $1$-solenoids is isomorphic to the finite product of
properly arranged $1$-solenoids and the isomorphism will yield equivalences
between the two corresponding families of linear flows [\textbf{Lemma \ref
{equi}}], so we need only consider the families of flows on the finite
product of properly arranged $1$-solenoids for the purposes of
classification. And two linear flows $\Phi _{\overline{P}}^\omega $ and $
\Phi _{\overline{P}}^{\omega ^{\prime }}\,$on a product of properly arranged 
$1$-solenoids are equivalent if and only if there is an is an automorphism $
A_{\overline{P}}$ with $a\omega ^{\prime }=A\left( \omega \right) $ for some 
$a\in \mathbb{R-}\left\{ 0\right\} $ by the above and by \textbf{Theorem \ref
{Class}}.

\subsection{Example: Classification in Dimension 2}

We only consider properly arranged $P\,$ and $Q$. We classify the linear
flows on $\sum_{\left( P,Q\right) }\overset{def}{=}\mathfrak{i}\left(
\sum_P\times \sum_Q\right) $. There are then 3 cases:

\begin{enumerate}
\item  $P=Q$

\item  $P>Q$ [i.e., $P\geq Q$ but not $Q\geq P$]

\item  $P$ and $Q$ are not comparable.
\end{enumerate}

\begin{corollary}
\label{equal}All equivalences of linear flows on $\sum_{\left( P,P\right) }$
are generated by automorphisms of the form $A_{P\times P}=\left( 
\begin{array}{cc}
a & b \\ 
c & d
\end{array}
\right) _{P\times P}$ for which all entries of $A_{P\times P}$ and $
A_{P\times P}^{-1}$ are proper $P\rightarrow P$ multipliers. In particular,
all rational linear flows on $\sum_{\left( P,P\right) }$ are equivalent.
\end{corollary}

\textbf{Proof}: The first statement follows from \textbf{Theorem \ref{aut}}.
Choosing appropriate integers $a,b,c$ and $d$ we can equate any two rational
flows with an automorphism of the form $\left( 
\begin{array}{cc}
a & b \\ 
c & d
\end{array}
\right) _{P\times P}$ where $\det $ $\left( 
\begin{array}{cc}
a & b \\ 
c & d
\end{array}
\right) =\pm 1$ (see \cite{I}, p.36).\hfill$\square $

The special case $P=(1,1,...)$ corresponds to the torus and the classical
classification of linear flows on the torus: in this case the only proper $
P\rightarrow P$ multipliers are integers and $f\left( Aut\left( \sum_{\left(
P,P\right) }\right) \right) =GL\left( 2,\mathbb{Z}\right) $.

\begin{corollary}
\label{bigger}If $P>Q\,$ are properly arranged, then all
equivalences of irrational flows on $\sum_{\left( P,Q\right) }$ may be
induced by automorphisms of the form $\left( 
\begin{array}{cc}
a & 0 \\ 
c & d
\end{array}
\right) _{(P,Q)},$ where $c$ and $\frac c{ad}$ are proper $P\rightarrow Q$
multipliers and $a$ and $d$ are iso-multipliers.\hfill$\square $
\end{corollary}

\begin{corollary}
\label{incomp}If $P$ and $Q\,$ are not comparable, then all equivalences of
linear flows on $\sum_{\left( P,Q\right) }$ may be induced by automorphisms
of the form $\left( 
\begin{array}{cc}
a & 0 \\ 
0 & d
\end{array}
\right) _{(P,Q)}$ , where $a$ and $d$ are iso-multipliers.
\end{corollary}

\textbf{Proof}: By the hypothesis $b$ and $c$ must be $0$ , since this is
the only proper $Q\rightarrow P$ $\left( P\rightarrow Q\right) $ multiplier.
Any such $\left( 
\begin{array}{cc}
a & 0 \\ 
0 & d
\end{array}
\right) _{(P,Q)}$ as stated is an automorphism with inverse $\left( 
\begin{array}{cc}
\frac 1a & 0 \\ 
0 & \frac 1d
\end{array}
\right) _{(P,Q)}$.\hfill$\square $

A somewhat surprising example of this case is given by $P=\left(
2,5,...\right) $ and $Q=\left( 3,7,...\right) $, the sequences of the odd
and even indexed primes. Then $f\left( Aut\left( \sum_{\left( P,Q\right)
}\right) \right) =\{\left[ 
\begin{array}{cc}
\pm 1 & 0 \\ 
0 & \pm 1
\end{array}
\right] \}\,$ since neither $P\,$ nor $Q$ has a prime divisor, implying that
at most two distinct phase portraits of linear flows are represented by any
given equivalence class of $\mathcal{F}_{\left( P,Q\right) }$.

In general, \textbf{Theorem} \textbf{\ref{aut} }may be used to determine
whether a given matrix represents an element of $Aut\left( \sum_{\left(
P_1,...,P_n\right) }\right) $. When $P_1=\cdots =P_n$ or when $P_1,...,P_n\,$
are pairwise incomparable, the classification works out as in \textbf{
Corollaries \ref{equal} }and \textbf{\ref{incomp} }respectively. However,
there are many other possibilities besides these in the general case
corresponding to the one case \textbf{\ref{bigger}\thinspace } as above when 
$n=2$. For example, when $n=3$ there are 7 additional distinct
possibilities. When $P_1=\cdots =P_n$ $=\left( 1,1...\right) $, we obtain
the well known result $f\left( Aut\left( \sum_{\left( P_1,...,P_n\right)
}\right) \right) =GL\left( n,\mathbb{Z}\right) $ since the only proper $
P_j\rightarrow P_i$ multipliers are integers, giving the classification of
linear flows on $\mathbf{T}^n$.

\section{Appendix}

Now we find simple conditions on the character group of an $n-$solenoid that
determine when it is a product of $1-$solenoids. For examples of $2-$
solenoids which are not such products see \cite{KM} and \cite{GR}. Recall that a 
\underline{$\mathbb{Q-}$basis} for a subgroup $G\subset \left( \mathbb{R},+\right) 
$ is a set $B=\left\{ \beta _i\right\} _{i\in I}$ satisfying the property
that all $g\in G-\left\{ 0\right\} $ can be uniquely represented as a sum $g=
\frac{c_1}{d_1}\beta _{i_1}+\cdots +\frac{c_n}{d_n}\beta _{i_n}$, where for $
i=1,...,n$ $\frac{c_i}{d_i}\in \mathbb{Q}-\left\{ 0\right\} $ and $\gcd \left(
c_i,d_i\right) =1$. And we refer to the representation $g=\sum_{i\in I}\frac{
c_i}{d_i}\beta _i$ with $\frac{c_i}{d_i}=\frac 01$ for all $0$ terms as the 
\underline{canonical representation} of $g$. We shall need the following
definition.

\begin{definition}
We define a set of generators $\mathcal{S}$ for a group $G\subset \left( 
\mathbb{R},+\right) $ to be \underline{relatively prime} with respect to the $
\mathbb{Q-}$basis $\left\{ \beta _1,...,\beta _n\right\} $ if for each $\lambda
\in \mathcal{S}$ with canonical representation 
\begin{equation*}
\lambda \overset{can}{=}\frac{c_1}{d_1}\beta _1+\cdots +\frac{c_n}{d_n}\beta
_n
\end{equation*}
we have $\gcd \left( d_j,\frac d{d_j}\right) =1$ for $j=1,...,n$, where $
d=d_1\times \cdots \times d_n$.
\end{definition}

\begin{theorem}
The $n-$solenoid $\sum_{\overline{N}}$ is isomorphic with a product of $n$ $
1-$solenoids $\Leftrightarrow $ its character group $\widehat{\sum_{
\overline{N}}}$ is isomorphic with a countable subgroup of the reals $G$ (in
the discrete topology) which has a set of generators $\mathcal{S}$ $=\left\{
\lambda _i\overset{can}{=}\frac{c_1^i}{d_1^i}\beta _1+\cdots +\frac{c_n^i}{
d_n^i}\beta _n\right\} $ which is relatively prime with respect to a $\mathbb{Q-
}$basis $B=\left\{ \beta _1,...,\beta _n\right\} $ $\subset $ $G$.
\end{theorem}

\textbf{Proof}:($\Rightarrow $) Suppose the conditions of the theorem are
met. Then we let 
\begin{equation*}
\;\delta _j^0\overset{def}{=}d_j^0\overset{def}{=}1\left( \text{for }
j=1,...,n\right) \,\text{and recursively}
\end{equation*}
\begin{eqnarray*}
\text{ for }i &\geq &1\text{ and all }j=1,...,n\text{ we let }\delta _j^i
\overset{def}{=}\text{lcm}(d_j^i,\delta _j^{i-1})\text{ and }\Delta _j^i
\overset{def}{=}\frac{\text{ }\delta _j^i}{\delta _j^{i-1}} \\
&&\text{and }M_i\overset{def}{=}\left( 
\begin{array}{ccc}
\Delta _1^i & 0 & 0 \\ 
0 & \ddots & 0 \\ 
0 & 0 & \Delta _n^i
\end{array}
\right) \text{ and }\overline{M}\overset{def}{=}\left( M_1,M_2,...\right) \\
\text{and} &&\mathcal{G}\overset{def}{=}\left\{ \beta _1,...,\beta _n,\frac{
\beta _1}{\delta _1^1},...,\frac{\beta _n}{\delta _n^1},...,\frac{\beta _1}{
\delta _1^i},\cdots ,\frac{\beta _n}{\delta _n^i},...\right\} \\
\text{and} &&\mathfrak{M}\overset{def}{=}the\text{ }subgroup\text{ }of\text{ }
\left( \mathbb{R},+\right) \text{ }generated\text{ }by\text{ }\mathcal{G}.
\end{eqnarray*}
Since $\mathfrak{M}$ is a subgroup of $\left( \mathbb{R},+\right) $ and 
\begin{equation*}
\frac{\delta _j^i}{d_j^i}\in \mathbb{Z},\text{ we have }\frac 1{d_j^i}\beta _j=
\frac{\delta _j^i}{d_j^i}\frac 1{\delta _j^i}\beta _j\in \mathfrak{M}\text{ for
any }i,j\text{.}
\end{equation*}
And so we have $\lambda _i\subset \mathfrak{M}$ for all $i$. Thus, $G\subset 
\mathfrak{M}$ and to show $G=\mathfrak{M}$ it suffices to show $G\supset \mathcal{G}$
. Assume inductively on $i\in \mathbb{N}$ that 
\begin{equation*}
\frac 1{d_j^{i-1}}\beta _j\text{ and }\frac 1{\delta _j^{i-1}}\beta _j\in G
\text{ for }j=1,...,n\text{ }
\end{equation*}
(this clearly holds for $i=1$). Let $d^i\overset{def}{=}d_1^i\cdots d_n^i$.
For $k\neq j,$ $
\begin{array}{c}
\underline{c_k^id^i} \\ 
d_j^id_k^i
\end{array}
$ $\in \mathbb{Z\Rightarrow }$
\begin{equation*}
\frac{c_j^id^i/d_j^i}{d_j^i}\beta _j=\frac{c_j^id^i}{d_j^id_j^i}\beta _j=
\frac{d^i}{d_j^i}\lambda _i-\frac{c_1^id^i}{d_j^id_j^1}\beta _1-\cdots -
\frac{c_{j-1}^id^i}{d_j^id_j^{j-1}}\beta _{j-1}-\frac{c_{j+1}^id^i}{
d_j^id_j^{j+1}}\beta _{j+1}-\cdots -\frac{c_n^id^i}{d_j^id_j^n}\beta _n\in G
\text{.}
\end{equation*}
Since $\mathcal{G}$ is relatively prime with respect to $B$, $\gcd ( 
\begin{array}{c}
\underline{c_j^id^i} \\ 
d_j^i
\end{array}
,d_j^i)=1$ and so there are integers $\mu $ and $\nu $ such that $\mu 
\begin{array}{c}
\underline{c_j^id^i} \\ 
d_j^i
\end{array}
+\nu d_j^i=1\Rightarrow $
\begin{equation*}
\frac 1{d_j^i}\beta _j=\frac{\mu c_j^id^i/d_j^i+\nu d_j^i}{d_j^i}\beta
_j=\mu \frac{c_j^id^i/d_j^i}{d_j^i}\beta _j+\nu \beta _j\in G\text{.}
\end{equation*}
Then we also have $\gcd ( 
\begin{array}{c}
\underline{\delta _j^i} \\ 
d_j^i
\end{array}
, 
\begin{array}{c}
\underline{\delta _j^i} \\ 
\delta _j^{i-1}
\end{array}
)=1$ and so there are integers $r$ and $s$ so that $r 
\begin{array}{c}
\underline{\delta _j^i} \\ 
d_j^i
\end{array}
+s 
\begin{array}{c}
\underline{\delta _j^i} \\ 
\delta _j^{i-1}
\end{array}
=1$. Hence, 
\begin{equation*}
\frac 1{\delta _j^i}\beta _j=\frac{r 
\begin{array}{c}
\underline{\delta _j^i} \\ 
d_j^i
\end{array}
+s 
\begin{array}{c}
\underline{\delta _j^i} \\ 
\delta _j^{i-1}
\end{array}
}{\delta _j^i}\beta _j=r\frac 1{d_j^i}\beta _j+s\frac 1{\delta
_j^{i-1}}\beta _j\in G,
\end{equation*}
completing the inductive step. Thus, $G=\mathfrak{M}$ and since $\mathfrak{M}$ is
isomorphic with the direct limit of 
\begin{equation*}
\mathbb{Z}^n\overset{M_1}{\longrightarrow }\mathbb{Z}^n\overset{M_2}{
\longrightarrow }\mathbb{Z}^n\overset{M_3}{\longrightarrow }\cdots ,
\end{equation*}
we have by Pontryagin duality that $\sum_{\overline{N}}\cong \sum_{\overline{
M}}$, which in turn is topologically isomorphic with a product of $n$
1-solenoids [\textbf{Lemma \ref{prod}}]. The other direction is clear since
duality respects finite products.\hfill$\square $

\end{document}